\documentclass[10pt,a4paper]{amsart}

\usepackage[centertags]{amsmath}

\usepackage{latexsym,amssymb}
\usepackage[ansinew]{inputenc}
\usepackage[T1]{fontenc}
\usepackage{times,mathptmx}
\usepackage{color, geometry}
\usepackage{hyperref}

\usepackage{amsfonts, times, mathptmx}
\usepackage{amssymb}
\usepackage{amsthm}
\usepackage{newlfont}

\usepackage{cite}

\newcommand {\real} {\mathbb{R}}

\theoremstyle{plain}
\newtheorem{thm}{Theorem} [section]
\newtheorem{cor}[thm]{Corollary}
\newtheorem{lem}[thm]{Lemma}
\newtheorem{prop}[thm]{Proposition}

\theoremstyle{definition}

\newcommand{\R}{\mathbb{R}}
\newcommand{\N}{\mathbb{N}}

\newcommand{\cC}{{\mathcal C}}

\newcommand{\cH}{{\mathcal H}}

\newcommand{\cK}{{\mathcal K}}      
\newcommand{\cL}{{\mathcal L}}   
\newcommand{\cM}{{\mathcal M}}   
\newcommand{\cN}{{\mathcal N}}

\newcommand{\cU}{{\mathcal U}}   
   
\newcommand{\cW}{{\mathcal W}}

\newcommand{\weak}{\rightharpoonup}

\newcommand{\eps}{\varepsilon}      
\renewcommand{\phi}{\varphi}

\DeclareMathOperator{\dist}{{dist}}
\DeclareMathOperator{\ext}{{ext}}
\DeclareMathOperator{\osc}{{osc}}
\DeclareMathOperator{\sgn}{{sgn}}

\begin {document}

\title[Least energy
nodal solutions to sublinear Neumann problems]{Existence, unique continuation and symmetry of least energy
nodal solutions to sublinear Neumann problems}

\author {Enea Parini, Tobias Weth}

\date{\today}

\address{Aix Marseille Universit\'{e}, CNRS, Centrale Marseille, I2M, UMR 7373, 13453 Marseille, France}

\email{enea.parini@univ-amu.fr}

\address{Institut f\"ur Mathematik, Goethe-Universit\"at Frankfurt, D-60054 Frankfurt am Main, Germany}

\email{weth@math.uni-frankfurt.de}



\begin {abstract}
We consider the sublinear problem
\begin {equation*}
\left\{\begin{array}{r c l c} -\Delta u  & = &|u|^{q-2}u & \textrm{in }\Omega, \\  u_\nu & = & 0 & \textrm{on
}\partial\Omega,\end{array}\right.
\end {equation*}
where $\Omega \subset \real^N$ is a bounded domain, and $1 \leq q <
2$. For $q=1$, $|u|^{q-2}u$ will be identified with $\sgn(u)$. We
establish a variational principle for least energy nodal solutions,
and we investigate their qualitative properties. In particular, we
show that they satisfy a unique continuation property (their zero set
is Lebesgue-negligible). Moreover, if $\Omega$ is radial, then least
energy nodal solutions are foliated Schwarz symmetric, and they are
nonradial in case $\Omega$ is a ball. The case $q=1$ requires special
attention since the formally associated energy functional is not
differentiable, and many arguments have to be adjusted.   
\end{abstract}

\maketitle

\section{Introduction}
\label{sec:neum-bound-cond}

Let $\Omega\subset\real^N$ a bounded open domain with Lipschitz boundary, and let $1\le q <2$. We are concerned with
the sublinear Neumann boundary value problem 
\begin {equation}
\label{neumann-problem}
\left\{\begin{array}{r c l c} -\Delta u  & = &|u|^{q-2}u & \textrm{in }\Omega, \\  u_\nu & = & 0 & \textrm{on
}\partial\Omega.\end{array}\right.
\end {equation} 
Here $u_\nu$ is the outer normal derivative of $u$ at the boundary $\partial
\Omega$, and the term $|u|^{q-2}u$ will be identified by $\sgn(u)$
in case $q=1$ in the following. For $q>1$, problem
\eqref{neumann-problem} arises e.g. in the study of the Neumann problem
for the (sign changing) porous medium equation. To see this, we set $v=|u|^{\frac{1}{m}-1}u$ with $m=\frac{1}{q-1}
> 1$ and note that \eqref{neumann-problem} may equivalently be written as  
\begin {equation}
\label{porous-problem}
\left\{\begin{array}{r c l c} -\Delta (|v|^{m-1}v)  & = &v & \textrm{in }\Omega, \\  v_\nu & = & 0 & \textrm{on
}\partial\Omega.\end{array}\right.
\end {equation} 
As a consequence, the function $w(x,t)=
[(m-1)t]^{-\frac{1}{m-1}} v(x)$ is a solution of the problem 
\begin {equation}
\label{neumann-problem-pme}
\left\{\begin{array}{r c l c} w_t -\Delta |w|^{m-1}w  & = &0 & \textrm{in }
    \Omega \times (0,
    \infty), \\  w_\nu & = & 0 & \textrm{on
}\partial\Omega \times (0,\infty).
\end{array}\right.
\end {equation} 
For more information on this relationship and a detailed discussion of
the (sign changing) porous medium equation, we refer the reader to
\cite[Chapter 4]{vazquez} and the references therein.\\
In the case
$q=1$, one may regard \eqref{neumann-problem} as a model problem within the
class of general elliptic boundary value problems with piecewise constant (and
therefore discontinuous) nonlinearities. Such problems
appear e.g. in the study of equilibria of reaction diffusion equations
with discontinuous reaction terms, see e.g. \cite{arrieta, chang:80,hilhorst.rodrigues:00,rauch}.\\
Integrating the equation in \eqref{neumann-problem} over
$\Omega$, we see that $\int_{\Omega} |u|^{q-2}u=0$ for every solution
of \eqref{neumann-problem}, hence every nontrivial solution is sign
changing. Let us consider the functional 
$$
\varphi: W^{1,2}(\Omega) \to \R, \qquad \varphi(u)= \frac{1}{2} \int_\Omega |\nabla u|^2\,dx -
\frac{1}{q}\int_\Omega |u|^q\,dx.  
$$ 
If $1<q<2$, $\varphi$ is of class $C^1$, and critical points of $\phi$
are precisely the weak solutions of \eqref{neumann-problem}. Moreover,
since the nonlinearity in \eqref{neumann-problem} is Hölder continuous,
weak solutions $u$ of \eqref{neumann-problem} are in
$C^{2,\alpha}_{loc}(\Omega)$ by elliptic regularity, and the
restriction of $u$ to the open set $\{u \not = 0\}$ is of class
$C^\infty$. If $q=1$, then $\phi$ fails to be differentiable and weak
solutions of \eqref{neumann-problem} are in general not of class
$C^2$, but they are still strong solutions contained in
$W^{2,p}_{loc}(\Omega)$ for every $p<\infty$ and thus contained in
$C^{1,\alpha}_{loc}(\Omega)$ for every $\alpha \in (0,1)$.

The purpose of this paper is to derive the existence of 
solutions of \eqref{neumann-problem} with minimal energy and to
characterize these solutions both
variationally and in terms of their qualitative properties. We first
consider the case $1<q<2$. In order to obtain least energy nodal solutions, we minimize the functional $\varphi$ on the set
\begin{equation}
  \label{eq:27}
\cN:= \Bigl \{u \in W^{1,2}(\Omega)\::\: \int_{\Omega}
|u|^{q-2}u\,dx=0 \Bigr\}
\subset W^{1,2}(\Omega).
\end{equation}
We shall see that minimizers of
$\phi|_\cN$ solve \eqref{neumann-problem}, so these minimizers are precisely the least energy nodal solutions of \eqref{neumann-problem}. Note that this property
does not follow from the Lagrange multiplier rules since $\cN$ is not
a $C^1$-manifold if $q<2$. Our main result for the case $1<q<2$ is the following.

\begin{thm}
\label{sec:neum-bound-cond-1}
Suppose that $1<q<2$, and let $m:= \inf \limits_{u \in \cN}\phi(u).$ 
Then we have:
\begin{itemize}
\item[(i)] $\varphi$ attains the value $m<0$ on $\cN$.
\item[(ii)] Every minimizer $u$ of $\varphi$ on $\cN$ is a sign changing solution of
\eqref{neumann-problem} such that $u^{-1}(0) \subset \Omega$ 
has vanishing Lebesgue measure. 
\item[(iii)] If $\Omega$ is a bounded radial domain, then every minimizer $u$
  of $\varphi$ on $\cN$ is foliated Schwarz symmetric.
\item[(iv)] If $\Omega=B_1(0)$ is the unit ball, then every minimizer $u$
  of $\varphi$ on $\cN$ is a nonradial function.
\end{itemize}
\end{thm}

We add some comments on these results.  The property $(i)$ follows by
standard arguments based on the weak lower continuity of the Dirichlet
integral $u \mapsto \int_{\Omega} |\nabla u|^2\,dx$. To show that
every minimizer of $\phi$ on $\cN$ is a solution of
\eqref{neumann-problem}, we use a saddle point characterization of $\cN$
(see Lemma~\ref{sec:neum-bound-cond-3} below). The most difficult part
is the unique continuation property of minimizers of $\phi$ on $\cN$,
i.e., the fact that their zero sets have vanishing Lebesgue
measure. Note that, due to the fact that the nonlinearity $u \mapsto
|u|^{q-2}u$ is not locally Lipschitz, the linear theory on unique
contination does not apply. Moreover, as can be seen from very simple
ODE examples already, nontrivial solutions of semilinear equations of
the type $-\Delta u = f(u)$ with non-Lipschitz $f$ may have very large
zero sets. It is an interesting open problem whether every nontrivial
solution of \eqref{neumann-problem} has the unique continuation property; we
conjecture that this is true. The proof of (iii) is again quite short
and essentially follows the arguments in
\cite{bartsch.weth.willem:05}. 
In contrast, the nonradiality property for least
energy nodal solutions stated in (iv) is not immediate. The idea is to use properties of
directional derivatives of $u$. For problems with
$C^1$-nonlinearities, nonradiality properties have
successfully been derived via
directional derivatives in the case of Dirichlet problems 
\cite{aftalion.pacella:04} and Neumann problems \cite{girao.weth:06},
while the methods in these papers differ significantly due to the impact of
the boundary conditions. A particular difficulty of the present problem is to analyze for which solutions the problem  \eqref{neumann-problem} has a meaningful linearization, see Proposition \ref{sec:neum-bound-cond-6} for a
first result on this question. 

Let us now consider the case $q=1$. In this case, the functional
$\varphi$ is not differentiable, so that the techniques used when
$1<q<2$ can not be applied. Moreover, the saddle point
characterization in Lemma~\ref{sec:neum-bound-cond-3} fails in the
case $q=1$, i.e., for the set $\cN=\{u \in W^{1,2}(\Omega)\::\:
\int_{\Omega}\sgn(u)\,dx = 0\}$. Nevertheless, we derive the same
conclusions as in Theorem~\ref{sec:neum-bound-cond-1} by adjusting the
variational principle. More precisely, we consider minimizers of
the restriction of $\varphi$ to the set
$$
 \cM := \bigl\{ u \in W^{1,2}(\Omega)\,:\,  \bigl| |\{u>0\}|-|\{u<0\}|
 \bigr| \le  |\{u=0\}| \bigr\}.
$$
Note that $\cM$ is strictly larger than $\cN$. Our main results for the case $q=1$
are collected in the following Theorem.

\begin{thm}
\label{sec:neum-bound-cond-8}
Suppose that $q=1$ in \eqref{neumann-problem} and the definition of
$\phi$, and let $m:= \inf \limits_{u \in \cM}\phi(u)$. Then we have:
\begin{itemize}
\item[(i)] $\varphi$ attains the value $m<0$ on $\cM$.
\item[(ii)] Every minimizer $u$ of $\varphi$ on $\cM$ is a sign changing solution of
\eqref{neumann-problem} such that $u^{-1}(0) \subset \Omega$ 
has vanishing Lebesgue measure. 
\item[(iii)] If $\Omega$ is a bounded radial domain, then every minimizer $u$
  of $\varphi$ on $\cM$ is foliated Schwarz symmetric.
\item[(iv)] If $\Omega=B_1(0)$ is the unit ball, then every minimizer $u$
  of $\varphi$ on $\cM$ is a nonradial function.
\end{itemize}
\end{thm}

The general strategy for the
proofs of (i)-(iii) is the same as in
Theorem~\ref{sec:neum-bound-cond-1}, but the details are quite
different due to the geometry of $\cM$ and the fact that $\phi$ fails
to be differentiable in the case $q=1$. The proof of (iv) is
completely different, since there seems to be no way to use
directional derivatives to prove nonradiality. Instead, our proof of
(iv) is based on inequalities comparing the value $m$ with
the least energy of {\em radial} nodal solutions of
\ref{neumann-problem} (in the case $q=1$). In fact, the latter value
can be computed explicitly once we have shown that least energy radial nodal solutions are strictly monotone in the radial variable and therefore have exactly two
nodal domains. We then compare this value with upper estimates for the value $m$
obtained by using the test functions $v(x)=x_1$ (if $n=2$) or $v(x)=\frac{x_1}{|x|}$ (if $n \geq 3$).

The paper is organized as follows. After proving some fundamental
properties of the functional $\varphi$ in the case $1<q<2$ (Section
\ref{sec:vari-fram-case}), we show that least energy nodal solutions
satisfy a unique continuation property (Section
\ref{sec:uniq-cont-prop}), and we then deal with symmetry results in
radially symmetric domains (Section \ref{sec:symmetry}). In
particular, Theorem~\ref{sec:neum-bound-cond-1} will readily follow
from Lemma~\ref{sec:neum-bound-cond-1-1} and 
Theorems~\ref{sec:neum-bound-cond-4},~\ref{sec:neum-bound-cond-5} and~\ref{sec:symmetry-results} below. In Section
\ref{sec:caseq1}, we turn to the case $q=1$ and prove Theorem~\ref{sec:neum-bound-cond-8}.

Finally, we mention that it is not straightforward to obtain
similar results for the Dirichlet problem corresponding
to \eqref{neumann-problem}. Indeed, least energy solutions of the
Dirichlet problem might have different variational characterizations
on different domains, so the situation is more complicated than in the
Neumann case. The Dirichlet problem is considered in a paper in preparation by the authors.\\  

{\bf Acknowledgement.} Part of this paper was written during several visits of E.P. to the Goethe-Universit\"{a}t in Frankfurt. He would like to thank the institution for its kind hospitality.

\section{The variational framework in the case $q>1$.}
\label{sec:vari-fram-case}
For fixed $q \in [1,2)$ and $u \in W^{1,2}(\Omega)$ with $\nabla u \not \equiv 0$, we put 
\begin{equation}
\label{eq:2}
 t^*(u)= \left(\frac{ \int_\Omega |u|^q}{\int_\Omega |\nabla
u|^2}\right)^{\frac{1}{2-q}} \quad \in \; (0,\infty).
\end{equation}
It is easy to see that $t^*(u)$ is the unique minimizer of the function 
\[
(0,\infty) \to \R,\qquad  t \mapsto \varphi(tu).
\]
Suppose that $1<q<2$ from now on, and consider the set $\cN$ defined
in (\ref{eq:27}). For any $u \in \cN \setminus \{0\}$, the value $t^*(u)$ is well defined and satisfies
$t^*(u)u \in \cN$ and $\phi(t^*(u)u)<0$. This in particular implies
that the infimum $m$ of $\phi|_{\cN}$ is negative. The next lemma
highlights the saddle point structure given by $\phi$ and the set $\cN$.

\pagebreak

\begin{lem}$ $
\label{sec:neum-bound-cond-3}
\begin{itemize}
\item[(i)] For every $u \in L^1(\Omega)$ there exists precisely one $c=c(u) \in \R$ such that
$\int_{\Omega}|u+c|^{q-2}(u+c)\,dx = 0$. Moreover, the map $L^1(\Omega) \to
\R$, $u \mapsto c(u)$ is continuous.
\item[(ii)] If $u \in W^{1,2}(\Omega)$, then 
\begin{equation}
  \label{eq:1}
\frac{\partial}{\partial c} \varphi(u+c)>0 \quad \text{for $c<c(u)$}
\qquad \text{and}\qquad \frac{\partial}{\partial c} \varphi(u+c)<0 \quad \text{for $c>c(u)$.}  
\end{equation}
 In particular,
\[ \varphi(u+c(u)) = \max_{c \in \real} \varphi(u+c).\]
\end{itemize}
\end{lem}

\begin{proof}
(i) Let $u \in L^1(\Omega)$. Since the map $\R \to \R,\; t \mapsto |t|^{q-2}t$
is strictly increasing, there exists at most one $c=c(u)$ such that
$\int_{\Omega}|u+c|^{q-2}(u+c)\,dx = 0$. Moreover, since 
$$
\lim_{c \to - \infty} \int_{\Omega}|u+c|^{q-2}(u+c)\,dx= -\infty\qquad
\text{and}\qquad \lim_{c \to \infty}
\int_{\Omega}|u+c|^{q-2}(u+c)\,dx= \infty,
$$
there exists precisely one such $c=c(u)$. Now consider $u, u_n  \in
L^1(\Omega)$, $n \in \N$ such that $u_n \to u$ in $L^1(\Omega)$. We
first show that $c(u_n)$ remains bounded as $n \to \infty$. Suppose by
contradiction that, after passing to a subsequence, $c(u_n) \to + \infty$ as $n \to
\infty$. Passing again to a subsequence, we may assume that $c(u_n)>0$
for all $n$ and $u_n \to u$ pointwise a.e. in $\Omega$. Moreover,
by~\cite[Lemma A.1]{willem} we may assume that there exists $\tilde u
\in L^1(\Omega)$ with $|u_n| \le \tilde u$ a.e. in $\Omega$ for all $n \in \N$. Since $-\tilde u \le u_n +c(u_n)$ a.e. in $\Omega$ for all
  $n \in \N$, we also have
$$
-|\tilde u|^{q-2}\tilde u \le
|u_n+c(u_n)|^{q-2}(u_n+c(u_n)) \qquad \text{a.e. in $\Omega$ for all
  $n \in \N$.} 
$$
Hence, since $|u_n+c(u_n)|^{q-2}(u_n+c(u_n)) \to \infty$ pointwise a.e. in
$\Omega$, Fatou's Lemma implies that 
$$
\int_{\Omega}|u_n+c(u_n)|^{q-2}(u_n+c(u_n)) \to \infty \qquad \text{as
    $n \to \infty$,}
$$
which contradicts the definition of the map $c$. In the same way, we
obtain a contradiction when assuming that $c(u_n) \to -\infty$ for a subsequence. Consequently, $c(u_n)$
remains bounded as $n \to \infty$. We now argue by contradiction,
supposing that $c(u_n) \not \to c(u)$ as
$n \to \infty$. Then we may pass to a subsequence such that $c(u_n)
\to c \not = c(u)$ as $n \to \infty$. Since the map 
$$
L^1(\Omega) \to \R,\qquad u \mapsto \int_{\Omega}|u|^{q-2}u
$$
is continuous and $u_n + c(u_n) \to u + c$ in $L^1(\Omega)$, we have that 
$$
\int_{\Omega}|u+c|^{q-2}(u+c)\,dx = 0.
$$
By the uniqueness property noted above, we then deduce that $c =
c(u)$, a contradiction. We thus conclude that $c(u_n) \to c(u)$ as $n
\to \infty$, and this shows the continuity of the map $c: L^1(\Omega)
\to \R$.\\ 
(ii) We have
$$
\frac{\partial}{\partial c} \varphi(u+c)= -
\int_{\Omega}|u+c|^{q-2}(u+c)\,dx \;
\left \{
  \begin{aligned}
  &>0 \qquad \text{for $c<c(u)$};\\
  &<0 \qquad \text{for $c>c(u)$},
  \end{aligned}
\right.
$$
as claimed. 
\end{proof}

\begin{lem}
\label{sec:neum-bound-cond-1-1}
The functional $\varphi$ attains the value $m<0$ on $\cN$. Moreover,
every minimizer $u$ of $\varphi$ on $\cN$ is a sign changing solution of
\eqref{neumann-problem}
\end{lem}

\begin{proof}
We first note that, as a consequence of
Lemma~\ref{sec:neum-bound-cond-3}(ii), we have  
$$
\|u\|_{L^q(\Omega)}^2 = \min_{c \in \R}\|u+c\|_{L^q(\Omega)}^2
\le |\Omega|^{2-q} \min_{c \in \R}\|u+c\|_{L^2(\Omega)}^2 \le
|\Omega|^{2-q}\mu^{-1}_2 \int_{\Omega}|\nabla u|^2\,dx \qquad \text{for $u
  \in \cN$,}
$$
where $\mu_2>0$ is the first nontrivial eigenvalue of the Neumann
Laplacian on $\Omega$. As a consequence, the functional $\varphi$ is coercive on $\cN$. Let
$(u_n)_n  \subset \cN$ be a minimizing sequence for $\varphi$. Then
$(u_n)$ is bounded, and we may pass to a subsequence such that $u_n
\weak u \in  W^{1,2}(\R^N)$. Then $u_n \to u$ in $L^q(\Omega)$,
$$
\int_{\Omega}|\nabla u|^2\,dx \le \liminf_{n \to \infty} |\nabla u_n|^2\,dx 
$$
and 
$$
\int_{\Omega}|u_n|^{q-2}u_n \,dx \to \int_{\Omega}|u|^{q-2}u \,dx
\qquad \text{as $n \to \infty$.}
$$
Consequently, we have $u \in \cN$, and $u$ satisfies $\varphi(u) \le
\liminf \limits_{n \to
  \infty}\varphi(u_n)$. Hence $u$ is a minimizer for $\varphi$ on $\cN$.\\
Next, we let $u \in \cN$ be an arbitrary minimizer for $\varphi$ on
$\cN$. We show that $u$ is a critical point of $\varphi$. Arguing again
by contradiction, we assume that there exists $v \in W^{1,2}(\Omega)$
such that $\varphi'(u)v < 0$. Since $\varphi$ is a $C^1$-functional on $W^
{1,2}(\Omega)$, there exists $\eps>0$ with
the following property: \medskip

{\em For every $w \in W^{1,2}(\Omega)$ with
  $\|w\|_{W^{1,2}(\Omega)}<\eps$ and every $t \in (0,\eps)$ we have
  \[\varphi(u+w+tv) \le \varphi( u +w)-\eps t.\]}

Since the map $c$ is continuous and $c(u)=0$ by definition of $c$, there exists $t \in (0,\eps)$ such
that $\|c(u+tv)\|_{W^{1,2}(\Omega)}<\eps$, and thus 
$$
\varphi(u+tv+c(u+tv))\leq \varphi(u+c(u+tv))-\eps t \leq  \varphi(u)-\eps t <\varphi(u).
$$
Since $u+tv+c(u+tv) \in \cN$, this contradicts the definition of
$m$. Finally, since $m<0$ by the remarks in the beginning of this
section, every minimizer $u \in \cN$ of $\phi$ is a nonzero function and
therefore sign changing by the definition of $\cN$.
\end{proof}

We close this section with a result on the existence of second
derivatives of $\phi$ which we will need in Section
\ref{sec:neum-bound-cond-5} below.

\begin{prop}
\label{sec:neum-bound-cond-6}
Let 
\begin{equation}
  \label{eq:15}
\cW:= \{v \in C^1(\overline \Omega)\::\: \text{$\nabla v(x) \not= 0$
  for every $x \in \Omega$ with $v(x)=0$}\}.
\end{equation}
Then $\cW \subset C^1(\overline \Omega)$ is an open subset (with respect
to the $C^1$-topology) having the following properties: 
\begin{itemize}
\item[(i)] If $\partial \Omega$ is of class $C^2$, then the restriction $\phi|_{\cW}$ is of class
$C^2$
with 
$$
\phi''(u)(v,w)= \int_{\Omega}\nabla v \nabla w\,dx - (q-1)
\int_{\Omega}|u|^{q-2}vw \qquad \text{for every $v,w \in C^1(\overline
  \Omega)$.}
$$
\item[(ii)] If $\partial \Omega$ is of class $C^{2,1}$, and $u \in \cW$ is a weak solution of \eqref{neumann-problem}, then
  $u \in W^{3,p}(\Omega)$ for $p \in (1,\frac{1}{2-q})$, and the partial derivatives 
$u_{x_i} \in W^{2,p}(\Omega)$ are strong solutions of the
problem 
\begin{equation}
  \label{eq:14}
-\Delta u_{x_i}= (q-1)|u|^{q-2}u_{x_i} \qquad \text{in $\Omega$.}  
\end{equation}
\end{itemize}
\end{prop}

\begin{proof}
It is easy to see that $\cW$ is open in $C^1(\overline \Omega)$. We
first 
show\\ 
{\em \underline{Claim 1:} If $s \in (0,1)$ and $1 \le p <
\frac{1}{s}$, then the map 
$$
\gamma_{s}: \cW \to L^p(\Omega), \qquad \gamma_{s}(u) \mapsto |u|^{-s}
$$
is well defined and continuous.}\\
To see this, let $\cK \subset \cW$ be a compact subset (with
respect to the $C^1$-norm). We claim that exists $\kappa>0$ such that 
\begin{equation}
  \label{eq:10}
|\{|u| \le \delta\}| \le \min \{\kappa \delta, |\Omega|\}  \qquad \text{for every
  $\delta>0$, $u \in \cK.$}
\end{equation}
In order to prove this estimate, we consider a bounded linear extension map $\ext: C^1(\overline
\Omega) \to C^1_b(\R^N)$, where $C^1_b(\R^N)$ denotes the Banach space of bounded $C^1$-functions
on $\R^N$ with bounded gradient. Such a map exists since $\partial
\Omega$ is of class $\cC^2$. We put $\cL:= \ext(\cK) \subset
C^1_b(\R^N)$ and $\Omega_b:= \{x \in \R^N \::\: \dist(x,\Omega) < b\}
\subset \R^N$ for $b>0$. Since $\cL$ is
compact, there exists a constant 
$b>0$ such that 
\begin{equation}
  \label{eq:19}
\max_{j=1,\dots,N}\Bigl|\frac{\partial v}{\partial x_j}(x)\Bigr| \ge 2b \qquad \text{for every $v \in \cL$ and $x
  \in \Omega_b$ with $|v(x)| \le b$.}
\end{equation}
We now fix $v_0 \in \cL$. Then there exists a positive integer
$d=d(v_0)$ and a finite number of cubes
$W_1,\dots,W_d$ of equal length $l>0$ such that
\begin{itemize} 
\item[(I)] every cube has the form $[x_1,x_1+l] \times \dots \times
  [x_N,x_N+l]$ with some $x \in \R^N$;
\item[(II)] $\overline \Omega \subset \bigcup \limits_{i=1}^d W_i \subset \Omega_b$; 
\item[(III)] $\underset{W_i}{\osc} \frac{\partial v_0}{\partial x_j} < \frac{b}{2}$ 
for $i=1,\dots,d$, $j=1,\dots,N$.
\end{itemize}
Moreover, there exists a neighborhood $\cU(v_0)$ of $v_0$ in $\cL$ such that
\begin{equation}
  \label{eq:21}
\underset{W_i}{\osc} \frac{\partial v}{\partial x_j} <b \qquad \text{ 
for $i=1,\dots,d$, $j=1,\dots,N$ and $v \in \cU(v_0)$.} 
\end{equation}
For every $v \in \cU(v_0)$, $i \in \{1,\dots,d\}$ and $\delta \in
(0,b)$ 
we then have 
\begin{equation}
  \label{eq:20}
|\{x \in W_i\::\: |v(x)| \le \delta\}| \le \frac{l^{N-1}}{b}\delta.
\end{equation}
Indeed, if there exists $x \in W_i$ with $|v(x)| \le \delta \le b$, then
$\Bigl|\frac{\partial v}{\partial x_j}\Bigr| \ge b$
on $W_i$ for some $j=j(i)$ by (\ref{eq:19}) and (\ref{eq:21}); in particular, $v$ is strictly monotone in the $j$-th coordinate
direction on $W_i$. Hence (\ref{eq:20}) easily follows by Fubini's theorem. 
As a consequence, we have the estimate 
$$
|\{x \in \overline \Omega\::\: |v(x)| \le \delta \}| \le 
\frac{d l^{N-1}}{b}\delta \qquad \text{for every $v \in \cU(v_0)$,
  $\delta \in (0,b)$.}
$$
Since $\cL$ is compact, it can be covered by finitely many
neighborhoods constructed as above, and hence there exists $d_*>0$
such that 
$$
|\{x \in \overline \Omega\::\: |v(x)| \le \delta \}| \le d_* \delta \qquad \text{for every $v \in \cL$, $\delta \in (0,b)$.}
$$
By the construction of $\cL$, (\ref{eq:10}) follows with $\kappa:=
\max \{d_*,\frac{|\Omega|}{b}\}$. As a consequence of (\ref{eq:10}),
we have  
\begin{align*}
\int_{\Omega}|u|^{-sp}\,dx = \int_0^\infty |\{|u|^{-sp} \ge \tau\}|\,d\tau =
\int_0^\infty |\{|u| \le \tau^{-\frac{1}{ps}}\}|\,d\tau \le \int_0^{\infty} \min
\{\kappa \tau^{-\frac{1}{ps}},|\Omega|\}\,d\tau < \infty,
\end{align*}
for every $u \in \cK$, since $\frac{1}{ps}>1$. In particular, the map $\gamma_s$ is well defined. To
see the continuity of $\gamma_s$, let $(u_n)_n
\subset \cW$ be a sequence such that $u_n \to u$ as $n \to \infty$
with respect to the $\cC^1$-norm. We then consider the compact set $\cK:=
\{u_n,u \::\: n \in \N\}$ and $\kappa>0$ such that (\ref{eq:10})
holds. For given $\eps>0$, we then fix $c>0$ sufficiently
small such that 
\begin{equation}
  \label{eq:9}
2^{p}  \int_{(2c)^{-ps}}^\infty
\min \{\kappa \tau^{-\frac{1}{ps}},|\Omega|\}\,d\tau<\eps.  
\end{equation}
By Lebesgue's theorem, it is easy to see
that 
\begin{equation}
  \label{eq:11}
\int_{|u| > c} (|u_n|^{-s}-|u|^{-s})^p \,dx \to 0 \qquad
\text{as $n \to \infty$.}
\end{equation}
Moreover, there exists $n_0 \in \N$ be such that $\{|u| \le c \} \subset
\{|u_n| \le 2c\}$ for $n \ge n_0$. Consequently,  
\begin{align*}
\int_{|u| \le c} \Bigl||u_n|^{-s}-|u|^{-s}\Bigr|^p \,dx  &\le 
2^{p-1} \int_{|u| \le c}\Bigl(|u_n|^{-ps}+ |u|^{-ps}\Bigr) \,dx \\ & \le 
2^{p-1} \Bigl(\int_{|u_n| \le 2 c} |u_n|^{-ps} \,dx + \int_{|u| \le c}|u|^{-ps}\,dx\Bigr)\\
&= 2^{p-1} \Bigl( \int_{(2c)^{-ps}}^\infty
|\{|u_n| \le \tau^{-\frac{1}{ps}}\}|\,d\tau + \int_{c^{-ps}}^\infty
|\{|u| \le \tau^{-\frac{1}{ps}}\}|\,d\tau\Bigr) \\
&\le 2^{p}  \int_{(2c)^{-ps}}^\infty
\min \{\kappa \tau^{-\frac{1}{ps}},|\Omega|\}\,d\tau < \eps \qquad \text{for $n \ge n_0$.}
\end{align*}
Combining this with (\ref{eq:11}), we conclude that 
$$
\limsup_{n \to \infty}\int_{\Omega} \Bigl||u_n|^{-s}-|u|^{-s}\Bigr|^p
\,dx \le \eps.
$$ 
Since $\eps>0$ was given arbitrarily, we conclude that 
$$
\|\gamma_{s}(u_n)-\gamma_{s}(u)\|_{L^p(\Omega)}^p = \int_{\Omega} \Bigl||u_n|^{-s}-|u|^{-s}\Bigr|^p
\,dx \to 0 \qquad \text{as $n \to \infty$.}
$$
Hence Claim 1 follows.\\
We now turn to the proof of (i). To show that $\phi|_{\cW}$ is of class $C^2$, it suffices to show that 
$$
\psi: \cW \to \R, \qquad \psi(u)= \frac{1}{q}\int_{\Omega}|u|^q\,dx
$$
is of class $C^2$ with 
\begin{equation}
  \label{eq:13}
\psi''(u)(v,w)= (q-1) \int_{\Omega}|u|^{q-2}vw \qquad \text{for every $v,w \in C^1(\overline
  \Omega)$.}
\end{equation}
By standard arguments, $\psi$ is of class $C^1$ with 
$$
\psi'(u)v= \int_{\Omega}|u|^{q-2}uv \qquad \text{for every $v \in C^1(\overline
  \Omega)$.}
$$
Let $u \in \cW$ and $v,w \in C^1(\overline \Omega)$ with
$\|v\|_{L^\infty(\Omega)}, \|w\|_{L^\infty(\Omega)}<1$. For $t \in \R
\setminus \{0\}$ we have 
$$
\frac{1}{t}\Bigl(\psi'(u+tw)v-\psi'(u)v\Bigr)=I_{t}+J_{t}
$$
with
$$
I_{t}= \int_{|u|> |t| }\frac{|u+tw|^{q-2}(u+tw)-|u|^{q-2}u}{t} v \,dx
,\quad  J_{t}= \frac{1}{t}\int_{|u|\le |t|}\frac{|u+tw|^{q-2}(u+tw) -
|u|^{q-2}u}{t}v\,dx.
$$ 
Note that, with $s:=2-q \in (0,1)$,  
\begin{align*}
I_{t}&= (q-1) \int_0^1 \int_{|u|> |t| }|u+\tau t w|^{-s}vw \,dx
\,d\tau\\
&= (q-1)\Bigl(\int_0^1 \int_{\Omega}\gamma_{s}(u+ \tau t w)vw\,dx \,d \tau -
 \int_0^1 \int_{|u|\le |t| }|u+\tau t w|^{-s}vw \,dx\,d\tau\Bigr)
\end{align*}
with 
$$
\int_0^1 \int_{\Omega}\gamma_{s}(u+ \tau t w)vw\,dx\,d\tau \to 
\int_{\Omega}\gamma_{s}(u)vw\,dx =
\int_{\Omega}|u|^{q-2}vw\,dx\qquad \text{as $t \to 0$}
$$
and, by applying H\"{o}lder's inequality with some $p \in (1,\frac{1}{s})$,  
$$
\Bigl|\int_0^1 \int_{|u|\le |t| }|u+\tau t w|^{-s} vw \,dx\,d\tau\Bigr| \le 
|\{|u| \le |t|\}|^{1-\frac{1}{p}} \int_0^1 \|\gamma_s(u+ \tau t w)\|_{L^p(\Omega)}\,d\tau
\to 0 \qquad \text{as $t \to 0$.}
$$
Moreover, by choosing $\kappa>0$ such that (\ref{eq:10}) holds for
$u$, we find that 
\begin{align*}
|J_{t}|&\le \frac{1}{|t|}\int_{|u|\le |t|} \Bigl(|u+tw|^{q-1}
+|u|^{q-1}\Bigr)\,dx\Bigr|=|t|^{q-2}\int_{|u|\le |t|}\Bigl|\frac{u}{t}+w\Bigr|^{q-1} +
\Bigl|\frac{u}{t}\Bigr|^{q-1}\,dx\Bigr|\\
&\le |t|^{q-2}(2^{q-1}+1)|\{|u| \le t\}| \le \kappa |t|^{q-1}(2^{q-1}+1) \to 0\qquad \text{as $t \to 0$.}
\end{align*}
Combining these estimates, we conclude that 
$$
\frac{1}{t}\Bigl(\psi'(u+tw)v-\psi'(u)v\Bigr) \to (q-1)
\int_{\Omega}|u|^{q-2}
vw \,dx  \qquad \text{as $t \to 0$.}
$$
Hence the second directional derivatives of $\phi$ exists at $u \in
\cW$ and satisfy (\ref{eq:13}). By Claim 1 above, it also follows that
the second derivatives depend continuously on $u \in \cW$, so that $\psi
\in C^2(\cW)$. The proof of (i) is thus complete.\\
To prove (ii), put $v:= |u|^{q-2}u$. Then for $f \in
\cC^\infty_c(\Omega)$ and $\eps>0$ sufficiently small we have, by the divergence theorem, 
$$
\Bigl|\int_{|u| \ge \eps}v \partial_{x_i} f\,dx - \int_{|u| \ge
  \eps}\partial_{x_i}v  f\,dx\Bigr| \le \int_{|u|=\eps}|v||f|\,d\sigma 
\to 0 \qquad \text{as $\eps \to 0$.}
$$
Hence 
$$
\int_{\Omega}v \partial_{x_i} f\,dx= \lim_{\eps \to
  0^+}\int_{|u| \ge \eps}v \partial_{x_i} f\,dx
= \lim_{\eps \to
  0^+}\int_{|u| \ge
  \eps}\partial_{x_i}v  f\,dx = (q-1) \int_{\Omega}|u|^{q-2}u_{x_i} f\,dx
$$
for $f \in \cC^\infty_c(\Omega)$, which shows that $v_{x_i}=
(q-1)|u|^{q-2}u_{x_i}$ in distributional sense. Considering $s= 2-q$ again, we deduce from Claim 1 above that $v_{x_i} \in
L^p(\Omega)$ for $p \in (1,\frac{1}{2-q})$.  Since moreover $-\Delta
u = v$ in $\Omega$, it follows from standard regularity theory that $u \in W^{3,p}(\Omega)$. Moreover, 
$$
-\Delta u_{x_i} = \partial_{x_i} (-\Delta u) = v_{x_i} \qquad \text{in
  $\Omega$}
$$
in strong sense for $i=1,\dots,N$, which shows (\ref{eq:14}).
\end{proof}

\section{The unique continuation property in the case $q>1$}
\label{sec:uniq-cont-prop}

In this section we still consider the case $1<q<2$, and we show that the set $u^{-1}(0) \subset \Omega$ has zero Lebesgue
measure for every minimizer of $\phi$ on $\cN$. For this we need some
preliminaries. We recall that, for a measurable subset $A \subset \R^N$, a point $x
\in \R^N$ is called {\em a point of density one for $A$} if 
$$
\lim_{r \to 0}\frac{|A \cap B_r(x)|}{|B_r(0)|} = 1.
$$
If $A \subset \R^N$ is measurable, then, by a classical result (see
e.g. \cite[p.45]{evans-gariepy}), a.e. $x \in A$ is a point of density one
for $A$. We
also need the following simple observation.

\begin{lem}
\label{sec:some-estim-relat-1}  
Let $\alpha>0$, and let $f: (0,\infty) \to \R$ be a nonnegative function with the following properties:
\begin{itemize}
\item[(i)] $f$ is bounded on $[\eps,\infty)$ for every $\eps>0$;
\item[(ii)] $f(r) = o(r^{-\alpha})$ as $r \to \infty$. 
\end{itemize}
Then for any $r>0$ there
exists $s>0$ with 
\begin{equation*}
f(s) \ge f(r)\qquad \text{and}\qquad   f(t) \le 2^{\alpha} f(s) \quad \text{for $t \in \left[\frac{s}{2},2s\right]$.}  
\end{equation*}
\end{lem}

\begin{proof}
We argue by contradiction. If the assertion was false, we would find
$r>0$ and a sequence $(s_n)_n \subset (0,\infty)$ such that $s_0=r$
and, for every $n \in \N$,
\begin{equation}
  \label{eq:4}
s_{n+1} \in \left[\frac{s_n}{2},2s_n\right] \qquad \text{and}\qquad f(s_{n+1})
> 2^{\alpha} f(s_n). 
\end{equation}
Without loss of generality, we may also assume that $f_0:= f(r)=
f(s_0)>0$ (otherwise we replace $s_n$ by $s_{n+1}$ for every $n$). From (\ref{eq:4}) we deduce 
\begin{equation}
  \label{eq:4-new}
s_n \ge 2^{-n} r \qquad \text{and}\qquad f(s_n) \ge 2^{n \alpha} f_0.
\end{equation}
Assumption (i) then implies that $s_n \to 0$ as $n
\to \infty$,  whereas (\ref{eq:4-new}) implies that $f(s_n)s_n^\alpha \ge r^{\alpha}f_0$ for all $n \in \N$. This contradicts the
 assumption $f(r)=o(r^{-\alpha})$ as $r \to 0$.
\end{proof}

\begin{prop}
\label{sec:some-estim-relat-2}
Let $u$ be a solution of 
\begin{equation}
  \label{eq:3}
-\Delta u = |u|^{q-2}u \qquad \text{in $\Omega$,}
\end{equation}
and suppose that $x_0 \in \Omega$ is a point of density one for the set
$u^{-1}(0) \subset \Omega$. Then $u(x) = o(|x|^{\frac{2}{2-q}})$
as $x \to x_0$.   
\end{prop}

Here we recall that, by elliptic regularity theory, a distributional
solutions $u$ of (\ref{eq:3}) contained in
$W^{1,2}_{loc}(\Omega)$ is in fact a classical solution in
$C^{2,\alpha}_{loc}(\Omega)$ for some $\alpha>0$. 

\begin{proof}[Proof of Proposition~\ref{sec:some-estim-relat-2}]
Without loss, we may assume that $x_0=0 \in \Omega$ and that $u$ is
bounded in $\Omega$ (otherwise we replace $\Omega$ by a compactly contained subdomain
containing $x_0$). We extend $u$ to all of $\R^N$ by setting $u \equiv
0$ on $\R^N \setminus \Omega$.
Since $0$ is a point of density one for the set $u^{-1}(0)$ and $u$ is continuous in $0$ and bounded on $\R^n$, 
the function
$$
f:(0,\infty) \to \R, \qquad f(r)= r^{-\frac{2}{2-q}}\sup \limits_{|x|=r}|u(x)|
$$
satisfies the assumptions of Lemma~\ref{sec:some-estim-relat-1} with $\alpha = \frac{2}{2-q}$. We first
show the following \medskip

{\em Claim:} The function $f$ is bounded on $(0,\infty)$. \medskip

Arguing by contradiction, we assume that there exists a sequence of radii $r_n
>0$, $n \in \N$ such that $f(r_n)
\to \infty$ as $n \to \infty$. By Lemma~\ref{sec:some-estim-relat-1},
we then find $s_n>0$,
$n \in \N$ such that, for all $n \in \N$, 
\begin{equation}
  \label{eq:5}
f(s_n) \ge f(r_n) \qquad \text{and} \qquad f(t) \le C f(s_n) \quad
\text{for $t \in \left[\frac{s_n}{2},2s_n\right],$}
\end{equation}
where $C:=2^{\frac{q}{2-q}}$. In particular, $f(s_n)
\to \infty$ as $n \to \infty$. By the definition of $f$, this implies
that $s_n \to 0$ as $n \to \infty$, so without loss we may assume that
$B_{2s_n}(0) \subset \Omega$ for all $n \in \N$. We now put $\Omega_0:= B_2(0)
\setminus \overline{B_{\frac{1}{2}}(0)}$ and define,  for $n \in \N$,
the functions $v_n:\Omega_0 \to \R$ as \[v_n(x)= \frac{s_n^{-\frac{2}{2-q}}
u(s_n x)}{f(s_n)}.\] The functions $v_n$ solve the equations
\begin{equation}
  \label{eq:6}
-\Delta v_n = f(s_n)^{q-2} |v_n|^{q-2}v_n \qquad \text{in $\Omega_0$},
\end{equation}
whereas
$$
|v_n(x)| = \frac{|x|^{\frac{2}{2-q}}(|x| s_n)^{-\frac{2}{2-q}}|u(s_n x)|}{f(s_n)} \le
\frac{2^{\frac{2}{2-q}} f(|x| s_n)}{f(s_n)} \le 2^{\frac{2}{2-q}} C  \qquad \text{for $x \in \Omega_0$, $n \in
\N$.}
$$
Moreover, there exists a sequence of points $x_n \in S^1:= \{y
\in \R^N\::\: |y|=1\}$ such that $|u(s_n x_n)|=\sup_{|x|=s_n}|u(x)|= s_n^{\frac{2}{2-q}}f(s_n)$ and hence $|v_n(x_n)| = 1$ for $n \in
\N$. Using (\ref{eq:6}), elliptic regularity theory and the fact that
$f(s_n)^{q-2} \to 0$ as $n \to \infty$, we may pass to a subsequence such that 
\begin{align*}
&x_n \to \bar x \in S^1,\\
&\text{$v_n \to v$ in $C^1_{loc}(\Omega_0)$,}
\end{align*}
where $v$ is a harmonic function in $\Omega_0$ such that $v(\bar x)= 1$. In particular, $v
\not \equiv 0$. On the other hand, since $0$ is a point of density one
for
the set 
$u^{-1}(0)$, the sets $A_n:= \{x \in B_2(0) \setminus B_{\frac{1}{2}}(0)\::\:
v_n(x) \not= 0\}$ satisfy  
$$
\frac{|A_n|}{|B_2(0)|} \le  \frac{|\{x \in B_{2s_n}\::\: u(x) \not=
    0\}|}{|B_{2s_n}(0)|} \to 0 \qquad \text{as $n \to \infty$}
$$
and thus $v \equiv 0$ a.e. in $B_2(0) \setminus
B_{\frac{1}{2}}(0)$. This is a contradiction, and thus the above claim
is true.\\
To finish the proof of the proposition, it thus remains to show that
$f(r) \to 0$ as $r \to 0$. Arguing again by contradiction, we assume
that there exists $\eps>0$ and a sequence of radii $r_n
>0$, $n \in \N$ such that $r_n \to 0$ as $n \to \infty$ and $f(r_n)
\ge \eps$ for all $n \in \N$. We may assume that $B_{2r_n}(0) \subset
\Omega$ for all $n \in \N$. We then consider $w_n: \Omega_0 \to \R$,
$w_n(x)= r_n^{-\frac{2}{2-q}}
u(r_n x)$ for $n \in \N$. It follows from the claim above that
the functions $w_n$ are uniformly bounded in $\Omega_0$.
Moreover, $w_n$ solves 
$$
-\Delta w_n = |w_n|^{q-2}w_n \qquad \text{in $B_2(0) \setminus
  B_{\frac{1}{2}}(0)$},
$$
and there exists a sequence of points $x_n \in S^1$, $n
\in \N$ such that $w(x_n)=f(r_n) \ge \eps$ for all $n \in \N$. Using elliptic regularity theory again, we may pass to a subsequence such that 
\begin{align*}
&x_n \to \bar x \in S^1,\\
&\text{$w_n \to w$ in $C^1_{loc}(\Omega_0)$,}
\end{align*}
where $w \in C^{2,\alpha}_{loc}(\Omega_0)$
is a weak solution of $-\Delta w = |w|^{q-2}w$ in $\Omega_0$ such that
$w(\bar x) \ge \eps$. In particular, $w \not \equiv 0$. 
However, by the same argument as above, the sets $A_n:= \{x \in B_2(0) \setminus B_{\frac{1}{2}}(0)\::\:
w_n(x) \not= 0\}$ satisfy $\lim \limits_{n \to \infty}|A_n|=0$ and
therefore $w \equiv 0$ a.e. in $B_2(0) \setminus
B_{\frac{1}{2}}(0)$. This yields a contradiction, and thus we conclude
that $f(r) \to 0$ as $r \to 0$, as required.
\end{proof}

Next, we consider the family of energy functionals
\begin{equation}
  \label{eq:24}
\phi_r: W^{1,2}(B_r(0)) \to \R,\qquad \phi_r(v)=
\frac{1}{2}\int_{B_r(0)}|\nabla v|^2\,dx -\frac{1}{q}
\int_{B_r(0)}|v|^q\,dx
\end{equation}
for $r>0$. We also consider the scaling map
\begin{equation}
  \label{eq:25}
T_r: W^{1,2}(B_1(0)) \to W^{1,2}(B_r(0)),\qquad [T_r v](x)= r^{\frac{2}{2-q}}v\left(\frac{x}{r}\right).
\end{equation}
We then have the following. 

\begin{prop}
\label{sec:unique-continuation-1}
Let $u \in W^{1,2}(\Omega)$ be a solution of (\ref{eq:3}) in $\Omega$,
and let $x_0 \in \Omega$ be a point of density one for $u^{-1}(0) \subset
\Omega$. Moreover, let $K \subset W^{1,2}_0(B_1(0))$ be a compact set
such that 
$$
c_K:= \sup_{v \in K}\phi_1(v)<0.
$$
Then there exists $r_0>0$ with the following property:\\
$B_{r_0}(x_0)$ is contained in $\Omega$, and for every $r \in (0,r_0)$ and every $v \in K$ we have $\phi(u+v_r)<\phi(u)$,
where $v_r \in W^{1,2}_0(\Omega)$ is defined by 
$$
v_r(x)= \left\{
  \begin{aligned}
   &[T_r v](x-x_0)&&\qquad x \in B_{r}(x_0);\\
   &0 &&\qquad x \in \Omega \setminus B_r(x_0).     
  \end{aligned}
\right.
$$
\end{prop}

\begin{proof}
Without loss, we may assume that $x_0=0 \in \Omega$. Let $v \in K$. We
then have 
\begin{align}
\nonumber
\phi(u+v_r)&=\phi(u)+\phi_r(v_r)+ \int_{B_r(0)}\nabla u \nabla v_r\,dx
- \frac{1}{q} \int_{B_r(0)}\Bigl(|u+v_r|^q-|u|^q-|v_r|^q\Bigr))\,dx\\
&\le \phi(u)+\phi_r(v_r)+ \int_{B_r(0)}\nabla u \nabla v_r\,dx +
\frac{1}{q} \int_{B_r(0)}\Bigl(|u|^q- q|v_r|^{q-2}v_r u\Bigr))\,dx   \label{eq:8}
\end{align}
where in the last step we used that 
$$
\int_{B_r(0)} \Bigl(|u+v_r|^q-|v_r|^q-q|v_r|^{q-2}v_r u\Bigr)\,dx \ge 0  
$$
by the convexity of the function $t \mapsto |t|^q$. We note that
$$
\int_{B_r(0)}\nabla u \nabla v_r\,dx = \int_{B_r(0)}|u|^{q-2}u v_r\,dx  
$$
and therefore
$$
\Bigl|\int_{B_r(0)}\nabla u \nabla v_r\,dx\Bigr| \le
\|u\|_{L^\infty(B_r(0))}^{q-1} \int_{B_r(0)}|v_r|\,dx =
o(r^{\frac{2(q-1)}{2-q}})r^{N+\frac{2}{2-q}} \int_{B_1(0)}|v|\,dx =
o(r^{N+ \frac{2q}{2-q}}) \quad \text{as $r \to 0$.}
$$
Moreover,
$$
\int_{B_r(0)}|u|^q\,dx = o(r^{\frac{2q}{2-q}})|B_r(0)|= o(r^{N+
  \frac{2q}{2-q}})
$$
and
$$
\Bigl|\int_{B_r(0)}|v_r|^{q-2}v_r u\,dx \Bigr| \le
\|u\|_{L^\infty(B_r(0))} \int_{B_r(0)}|v_r|^{q-1} \,dx \le
o(r^{\frac{2}{2-q}}) r^{N+\frac{2(q-1)}{2-q}}
\int_{B_1(0)}|v|^{q-1}\,dx = o(r^{N+ \frac{2q}{2-q}}) 
$$
as $r \to 0$. Finally, since $\nabla v_r(x)= r^{\frac{q}{2-q}}\nabla v(\frac{x}{r})$ for $x
\in B_r(0)$, we find
that   
\begin{align*}
\phi_r(v_r)&= \frac{r^{\frac{2q}{2-q}}}{2} \int_{B_r(0)}\Bigl|\nabla
v\left(\frac{x}{r}\right)\Bigr|^2\,dx
-\frac{r^{\frac{2q}{2-q}}}{q}\int_{B_r(0)}\Bigl|v\left(\frac{x}{r}\right)\Bigr|^q\,dx\\
&=r^{N+ \frac{2q}{2-q}}\Bigl(\frac{1}{2} \int_{B_1(0)}|\nabla
v(x)|^2\,dx
-\frac{1}{q}\int_{B_1(0)}|v(x)|^q\,dx\Bigr)\\
&=r^{N+ \frac{2q}{2-q}}\phi_1(v).
\end{align*}
Inserting these estimates in (\ref{eq:8}), we obtain 
$$
\phi(u+v_r)\le \phi(u)+ r^{N+ \frac{2q}{2-q}}\Bigl(\phi_1(v)+o(1)\Bigr) \le \phi(u)+ r^{N+ \frac{2q}{2-q}}\Bigl(c_K+o(1)\Bigr) 
$$
Since $K$ is compact, it is easy to see that these estimates are
uniform in $v \in K$. Since moreover $c_K<0$, the claim follows.
\end{proof}

\begin{thm}
\label{sec:neum-bound-cond-4}
Let $u$ be a
minimizer of $\phi$ on $\cN$. Then $u^{-1}(0)$ has vanishing Lebesgue
measure.   
\end{thm}

\begin{proof}
Suppose by contradiction that $|u^{-1}(0)|>0$. Then there exists a
point $x_0$ of density one for the set $u^{-1}(0)$. Without loss of generality, we can suppose $x_0=0$. We fix two arbitrary nonnegative nontrivial functions $v_1,v_2 \in
W^{1,2}_0(B_1(0))$ with disjoint support, and consider the path 
$$
\gamma:[0,1] \to W^{1,2}_0(B_1(0)),\qquad \gamma(s)= t^{*}((1-s) v_1 -
s v_2)\cdot ((1-s) v_1 -s v_2),
$$
where the function $t^*$ is defined in (\ref{eq:2}). It is clear that $\varphi(\gamma(t))<0$ for $t \in [0,1]$. Applying
Proposition~\ref{sec:unique-continuation-1} to the compact set $K:=
\gamma([0,1]) \subset W^{1,2}_0(B_1(0))$, we may fix $r>0$
sufficiently small such that $\phi(u+ v_r)<\phi(u)$ for every $v \in
K$. Since $K$ is connected and 
$$
\int_{\Omega}|u+v_r|^{q-2}(u+v_r)\,dx \;
\left \{
  \begin{aligned}
  &>0 \qquad \text{for $v=\gamma(0) \in K$,}\\
  &<0 \qquad \text{for $v=\gamma(1) \in K$,}  
  \end{aligned}
\right.
$$
there exists $v \in K$ such that $u+v_r \in \cN$. This however
contradicts the assumption that $u$ is a minimizer of $\phi$ in $\cN$.
\end{proof}

\begin{cor}
\label{sec:neum-bound-cond-2}
Let $u$ be a
minimizer of $\phi$ on $\cN$, and let $x_0 \in \Omega$ be a point with
$u(x_0)=0$. Then $u$ changes sign in every neighborhood of $x_0$.   
\end{cor}

\begin{proof}
Suppose by contradiction that there is a ball $B_r(x_0) \subset
\Omega$ such that $u$ does not change sign in $B_r(x_0)$. Without
loss, we may assume that $u \ge 0$ on $B_r(x_0)$. By
Theorem~\ref{sec:neum-bound-cond-4}, $u \not \equiv 0$ in
$B_r(x_0)$. Since $u$ solves \eqref{neumann-problem} and is therefore
superharmonic in $B_r(x_0)$, the strong maximum principle implies that
$u>0$ in $B_r(x_0)$, contrary to the assumption that $u(x_0)=0$.  
\end{proof}

\section{Symmetry results}
\label{sec:symmetry}

We add a result on minimizers of $\phi|_{\cN}$ in the case
where the underlying domain is radial, i.e., a ball or an annulus in
$\R^N$ centered at zero. 

\begin{thm}
\label{sec:neum-bound-cond-5}
Suppose that $\Omega \subset \R^N$ is a radial bounded domain. Then
every minimizer $u$ of $\phi$ on $\cN$ is foliated Schwarz
symmetric.  
\end{thm}

Here we recall that a function $u$ defined on a
radial domain is said to be {\em foliated Schwarz symmetric} if there is a unit vector $p\in \R^N$,
 $|p|=1$ such that $u(x)$ only depends on $r=|x|$ and $\theta
 =\arccos \left(\frac x{|x|}\cdot p\right)$ and $u$ is nonincreasing
 in $\theta$.

 \begin{proof}
Let $u \in \cN$ be a minimizer of $\phi|_{\cN}$, and pick $x_0 \in
\Omega \setminus \{0\}$ with $u(x_0)=\max \{u(x)\::\: |x|=|x_0|\}$.
We put $p:=\frac{x_0}{|x_0|}$, and we let $\cH_p$ denote the
family of all open halfspaces $H$ in
$\R^N$ such that $p \in H$ and $0 \in \partial H$. For $H \in \cH_p$ we consider the reflection $\sigma_H: \R^N \to \R^N$
with respect to the the hyperplane $\partial H$. We claim the following:
\begin{equation}
  \label{eq:7}
\text{For every $H \in \cH_p$, we have $u \ge u \circ \sigma_H$ on $H \cap \Omega$.} 
\end{equation}
To prove this, we fix $H \in \cH_p$ and recall a simple rearrangement,
namely the {\em polarization} of $u$ with respect to $H$ defined by
$$
u_H(x)=
\begin{cases}
 &\max \{u(x),u(\sigma_H(x))\},\qquad x\in \Omega \cap H\\
 &\min \{u(x),u(\sigma_H(x))\},
\qquad x\in \Omega \setminus H.
\end{cases}
$$
It is well known and fairly easy to prove (see e.g. \cite{willem2007}) that
$$
\int_{\Omega} |\nabla u_H|^2\:dx = \int_{\Omega} |\nabla
u|^2\:dx,\quad
\int_{\Omega} |u_H|^q\:dx = \int_{\Omega} |u|^q\:dx \quad
\text{and}\quad \int_{\Omega} |u_H|^{q-2}u_H\:dx = \int_{\Omega}
|u|^{q-2}u \:dx.
$$
Consequently, $u_H \in \cN$ and $\phi(u_H)=\phi(u)$, so that $u_H$ is
also a minimizer of $\phi$ on $\cN$. Hence, by Theorem~\ref{sec:neum-bound-cond-1}, both $u$ and $u_H$ are
solutions of \eqref{neumann-problem}. Therefore $w:=u_H-u$ is a
nonnegative function in $\Omega \cap H$ satisfying 
$$
-\Delta w= |u_H|^{q-2}u_H - |u|^{q-2}u \ge 0\qquad \text{in
    $H \cap \Omega$.}
$$
The strong maximum principle then implies that either $w \equiv 0$ or $w>0$ in $H \cap \Omega$. The
latter case is ruled out since $x_0 \in H \cap \Omega$ and
$w(x_0)=u_H(x_0)-u(x_0)=0$ by the choice of $x_0$. We therefore obtain
$w \equiv 0$, hence $u=u_H$ and (\ref{eq:7}) holds.\\
By continuity, it follows from (\ref{eq:7}) that $u$ is symmetric with respect to every
hyperplane containing $p$, so it is axially symmetric with respect
to the axis $p\R$. Hence $u(x)$ only depends on $r=|x|$ and $\theta
 =\arccos \left(\frac x{|x|}\cdot p\right)$. Moreover, it also
 follows from (\ref{eq:7}) that $u$ is nonincreasing
 in the polar angle $\theta$. We thus conclude that $u$ is foliated Schwarz symmetric.
\end{proof}

\begin{thm}
\label{sec:symmetry-results}
Let $\Omega \subset \real^N$ be the unit ball, and let $u \in \cN$ be a minimizer of $\phi|_\cN$. Then $u$ is not radially symmetric. 
\end{thm}
\begin{proof}
By Lemma~\ref{sec:neum-bound-cond-1-1}, $u$ is a solution of \eqref{neumann-problem}, so $u \in C^{2,\alpha}(\overline \Omega)$ by elliptic regularity. We suppose by contradiction
that $u$ is radially symmetric, and we write $u(r):=u(|x|)$ for
simplicity. Then $u$ solves
\begin{equation}
  \label{eq:17}
u'' + \frac{N-1}{r}u'+|u|^{q-2}u =0 \qquad \text{in $(0,1]$,}\qquad 
u'(0)=u'(1)=0,
\end{equation}
where the prime denotes the radial derivative. We first prove\\
\textit{\underline{Claim 1:}} $u(0) \not=0$, $u(1) \not=0$, and $u$ has
only finitely many zeros in $(0,1)$.\\
Indeed, with the transformation
$$
y(t)=u(r) \qquad \text{with}\quad t= 
\left \{
  \begin{aligned}
  r^{2-N},\qquad N \ge 3;\\
  1-\log r,\qquad N=2,  
  \end{aligned}
\right.
$$
(\ref{eq:17}) is transformed into the problem 
\begin{equation}
  \label{eq:18}
\ddot y +
p(t)|y|^{q-2}y=0,\quad t \in [1,\infty),\qquad \dot y(1)=0
\end{equation}
with 
$$
p(t)= 
\left\{
  \begin{aligned}
  &\frac{1}{(N-2)^2} t^{-\frac{2(N-1)}{N-2}}, &&\qquad N
\ge 3,\\
  &e^{-2(t-1)},&&\qquad N=2. 
\end{aligned}
\right.
$$    
Since $u \not \equiv 0$ we have $y \not \equiv 0$ in $[1,\infty)$. 
Hence $y(1) \not= 0$ by local uniqueness and continuability of
solutions to (\ref{eq:18}) (see e.g. \cite{wong1975}), and thus $u(1)
\not=0$. Moreover, by the non-oscillation criterion for (\ref{eq:18}) given in \cite[Theorem 6]{gollwitzer}, $y$ has only finitely many zeros in
$[1,\infty)$, so $u$ only has finitely many zeros in
$(0,1)$. Finally, since $u$ does not change sign in some neighborhood of $0$, the strong
maximum principle implies that $u(0) \neq 0$. Thus Claim 1 is
proved.\\
It now follows from Hopf's boundary lemma that $u \in \cW$, where $\cW$ is defined in
(\ref{eq:15}). Consequently, $u \in W^{3,p}(\Omega)$ for $p \in (1,\frac{1}{2-q})$ by
Proposition~\ref{sec:neum-bound-cond-6}, and $u_{x_1} \in
W^{2,p}(\Omega) \cap C^{1}(\overline \Omega)$ solves the linearized Dirichlet problem
\begin {equation}\label{problemu1}\left\{\begin{array}{r c l c} -\Delta u_{x_1}  & = & (q-1)|u|^{q-2}u_{x_1} & \text{in }\Omega, \\ u_{x_1} & = & 0 & \text{on
}\partial\Omega.\end{array}\right.
\end {equation}
The boundary condition follows from the fact that $\nabla u \equiv 0$
on $\partial \Omega$ since $u$ is radial and satisfies Neumann
boundary conditions. Let $\mathcal{H}$ be the hyperplane $\{x_1=0\}$. We first prove the following\\
\textit{\underline{Claim 2:}} If $w \in C^1(\overline \Omega)$ is antisymmetric
with respect to $\mathcal{H}$ and such that $\varphi''(u)(w,w) < 0$, then $\varphi(u+tw+c(u+tw)) < \varphi(u)$ for $t>0$ sufficiently small.\\
Indeed, by Proposition~\ref{sec:neum-bound-cond-6}(i) 
we have, for every $c \in \R$, the Taylor expansion 
\[ \varphi(u+c+tw) = \varphi(u+c) + t \varphi'(u+c)w + \frac{t^2}{2} \varphi''(u+c)(w,w) + o(t^2), \]
where the quantity $o(t^2)$ is locally uniform in
$c$. Since $u$ is radially symmetric and $w$ is antisymmetric with
respect to $\mathcal{H}$, we have 
$$
\varphi'(u+c)w= \int_{\Omega}\nabla u \nabla w\,dx -
\int_{\Omega}|u+c|^{q-2}(u+c)w\,dx =0.
$$
Hence there exist $M>0$ and $\delta >0$ such that
$$
\varphi(u+tw+c) \leq \varphi(u+c) - M t^2 \leq \varphi(u) - Mt^2
\qquad \text{for $|t|,|c| < \delta$.}
$$
Since $c(u+tw) \to c(u)=0$ as $t \to 0$ as a consequence of
Lemma~\ref{sec:neum-bound-cond-3}, we deduce that 
$$
\varphi(u+tw+c(u+tw)) \leq  \varphi(u) - Mt^2 < \varphi(u) \qquad \text{for $t>0$ sufficiently small.}
$$
Hence Claim 2 is proved. Next, we consider an arbitrary 
function $w \in C^1(\overline \Omega)$ which is antisymmetric with
respect to $\cH$. By Claim 2 and the
minimizing property of $u$, we have
$$ 
0 \leq \varphi''(u)(w+tu_{x_1},w+tu_{x_1}) = \varphi''(u)(w,w) + 2t
\varphi''(u)(u_{x_1},w) + t^2 \varphi''(u)(u_{x_1},u_{x_1}) \qquad \text{for every
  $t \in \R$}
$$
and also 
\[ \varphi''(u)(u_{x_1},u_{x_1}) = \int_\Omega |\nabla u_{x_1}|^2\,dx - (q-1)\int_\Omega |u|^{q-2}u_{x_1}^2\,dx = 0, \]
since $u_{x_1}$ is a solution of \eqref{problemu1}. These relations imply that 
\begin{equation}
  \label{eq:16}
0= \varphi''(u)(u_{x_1},w) = \int_\Omega \nabla u_{x_1} \nabla w\,dx - (q-1)\int_\Omega |u|^{q-2}u_{x_1} w\,dx = \int_{\partial \Omega} (u_{x_1})_\nu w\,d\sigma,
\end{equation}
where the last equality follows again from \eqref{problemu1}. Since
$(u_{x_1})_\nu \in C(\partial \Omega)$ is antisymmetric with respect to
$\cH$ and (\ref{eq:16}) holds for every $w \in C^1(\overline \Omega)$
which is antisymmetric with respect to $\cH$, we conclude that
$(u_{x_1})_\nu = 0$ on $\partial \Omega$, and in particular $u_{x_1,x_1}
(e_1)=0$. In the radial variable, we thus have 
$$
0 = u''(1)= -|u(1)|^{q-2}u(1) 
$$
by (\ref{eq:17}) and therefore $u(1)=0$, contrary to Claim 1. The
proof is finished.
\end{proof}

\section{The case $q=1$}
\label{sec:caseq1}

In this section we are concerned with the case $q=1$, i.e., with the boundary value problem  
\begin {equation}
\label{problem_q=1}
\left\{\begin{array}{r c l c} -\Delta u  & = &\text{sgn}(u),  & \textrm{in }\Omega \\ u_\nu & = & 0 & \textrm{on
}\partial\Omega.\end{array}
\right.
\end {equation}
We will suppose that the boundary of $\Omega$ is of class $C^{1,1}$. As already noted in the introduction, the variational framework of
Section~\ref{sec:vari-fram-case} does not extend in a straightforward
way to the case $q=1$. In particular, the functional 
$$
\phi:W^{1,2}(\Omega) \to \R, \qquad \phi(u)= \frac{1}{2} \int_{\Omega}
|\nabla u|^2\,dx -\int_{\Omega}|u|\,dx 
$$
is not differentiable. Moreover, while every (weak) solution of
\eqref{problem_q=1} is contained in the set 
$$
\{ u \in W^{1,2}(\Omega)\,|\,\int_\Omega
\text{sgn}(u)\,dx=0\},
$$
this set does not have the nice intersection property given by
Lemma~\ref{sec:neum-bound-cond-3}(i). Indeed, any function 
$u \in W^{1,2}(\Omega)$ satisfying 
$$
0< \Bigl| |\{u>0\}|-|\{u<0\}| \Bigr| < |\{u=0\}|
$$
has the property that 
$$
\int_\Omega
\text{sgn}(u+c)\,dx \not=0 \qquad \text{for every $c
\in \R$.} 
$$
As a consequence, many of the arguments used in the previous sections do
not apply in the case $q=1$. Instead, we will consider the larger set
$$
 \cM := \bigl\{ u \in W^{1,2}(\Omega)\,:\,  \bigl| |\{u>0\}|-|\{u<0\}|
 \bigr| \le  |\{u=0\}| \bigr\}.
$$
We collect useful properties of $\cM$. First, we may rewrite the
defining property for $u \in \cM$ as 
\begin{equation}
  \label{eq:29}
\int_{\Omega} \sgn_-(u)\,dx \le 0 \le \int_{\Omega}\sgn_+(u)\,dx,  
\end{equation}
where 
$$
\sgn_+(t):= 1_{t \ge 0}- 1_{t<0} \quad \text{and}\quad \sgn_-(t):= 1_{t
  > 0}- 1_{t\le 0} \qquad \text{for $t \in \R$.}
$$
We also point out that, in contrast to the definition in the case $q>1$, the
set $\cM$ also contains nonzero functions which do not change sign. We also need the following facts.

\begin{lem}
\label{sec:case-q=1}  
\begin{itemize}
\item[(i)] If $u \in \cM$, then 
$$
\int_{\Omega}|u+c|\,dx > \int_{\Omega}|u|\,dx \qquad \text{for every $c
\in \R \setminus \{0\}$.}
$$
\item[(ii)] For $u \in W^{1,2}(\Omega)$, we have $u \in \cM$ if and only if 
$$
\int_{\Omega} \sgn(u-c)\,dx \le 0 \le \int_{\Omega} \sgn(u+c)\,dx \quad \text{for every $c>0$.}
$$
\item[(iii)] If $u \in W^{1,2}(\Omega)$ and $\gamma: [0,1] \to
  W^{1,2}(\Omega) \cap L^\infty(\Omega)$ is a continuous curve such
  that 
$$
\int_{\Omega}\sgn(u+\gamma(0))\,dx>0 \qquad \text{and}\qquad 
\int_{\Omega}\sgn(u+\gamma(1))\,dx<0,
$$
then there exists $s_0 \in [0,1]$ with $u + \gamma(s_0) \in \cM$.
\item[(iv)] For every $u \in W^{1,2}(\Omega)$ there exists a unique $c(u) \in
  \R$ such that $u+c(u) \in \cM$. Moreover, the map $u \mapsto c(u)$
  is continuous. 
\item[(v)] If $u,v \in W^{1,2}(\Omega)$ satisfy $u \le v$, then $c(u)
  \ge c(v)$.  
\end{itemize}
\end{lem}

\begin{proof}
(i) Let $c<0$. If $u \le 0$, then obviously 
$\int_{\Omega}|u+c|\,dx > \int_{\Omega}|u|\,dx$. Suppose now that $u^+
\not =0$. We then claim that 
\begin{equation} 
\label{eq:32} |\{0<u<|c|\}|>0 .
\end{equation} 
Indeed, we have $v:= \min (u^+,|c|) \not \equiv
0$, since $v>0$ on the set where $u$ is positive. Moreover, if we suppose by contradiction
that $|\{0<u<|c|\}|=0$, then $\nabla v \equiv 0$ a.e. on $\Omega$ (see
e.g.  \cite[Lemma 7.7]{gt}), so that $v$ equals a positive constant a.e. in $\Omega$. This however implies that $u>0$ a.e. in $\Omega$, contrary to the assumption $u \in \cM$. Hence we conclude that (\ref{eq:32}) holds. As a consequence,
we estimate, for $c<0$,
\begin{align*}
\nonumber
\int_{\Omega}|u+c|\,dx&-\int_{\Omega}|u|\,dx = |c| |\{u \le 0\}| +
\int_{\{0<u<|c|\}}(|u+c|-u)\,dx+
\int_{\{u\ge |c|\}}(|u+c|-u)\,dx\\
&= |c| |\{u \le 0\}| +
\int_{\{0<u<|c|\}}(|c|-2u)\,dx-|c| |\{u\ge |c|\}|\\
&> |c| |\{u \le 0\}| -|c| |\{0<u<|c|\}|-|c| |\{u\ge |c|\}|\\
&= |c| \bigl(|\{u \le 0\}| -|\{u >
0\}|\bigr)\ge 0,
\end{align*}
as claimed. If $c>0$, a similar argument yields 
$$
\int_{\Omega}|u+c|\,dx-\int_{\Omega}|u|\,dx >0.\\ 
$$   
(ii) This simply follows from the fact that $u \in \cM$ is equivalent
to (\ref{eq:29}), whereas 
$$
\sgn(u+c) \ge \sgn_+(u)  \ge \sgn_-(u) \ge \sgn(u-c) \qquad \text{for
  every $c>0$}
$$
and 
$$
\sgn(u+c) \to \sgn_+(u),\; \sgn(u-c) \to \sgn_-(u)\qquad
\text{pointwise in $\Omega$ as $c \to 0^+$.}
$$
(iii) Consider
$$
s_0:= \sup \bigl \{s \in [0,1)\::\: \int_{\Omega}\sgn(u+\gamma(s))\,dx
> 0\bigr \},
$$
and let $w:= u+\gamma(s_0)$. We use (ii) to show that $w \in \cM$. Let $(s_n)_n \subset [0,s_0]$
be a sequence with $s_n \to s_0$ and 
\begin{equation}
  \label{eq:31}
 \int_{\Omega}\sgn(u+\gamma(s_n))\,dx >0 \qquad \text{for every $n \in
  \N$.}
\end{equation}
For given $c>0$, there exists $n \in \N$ with 
$$
\{u+\gamma(s_n)>0\} \subset \{w+c>0\} \quad \text{and}\quad 
\{w+c<0\} \subset \{u+\gamma(s_n)<0\} 
$$
and thus $\int_{\Omega}\sgn(w +c)\,dx>
0$ by (\ref{eq:31}). Now if $s_0=1$, the assumption implies that 
$$
0> \int_\Omega \sgn(w)\,dx \ge \int_\Omega \sgn(w-c)\,dx \qquad
\text{for all $c>0$}
$$
and thus $w \in \cM$ by (ii).  Suppose finally that $s_0<1$, and suppose by
contradiction that 
$$
\int_{\Omega}\sgn(w -c)\,dx> 0 \qquad \text{for some $c>0$.}
$$
By the continuity of $\gamma$, there exists $\eps>0$ such that $u
+\gamma(s) \ge w-c$ for $s \in [s_0,s_0+\eps)$ and therefore 
$$
\int_{\Omega}\sgn(u+ \gamma(s))\,dx> 0 \qquad \text{for $s \in [s_0,s_o+\eps)$.}
$$
This contradicts the definition of $s_0$. Hence 
$$
\int_{\Omega}\sgn(w -c)\,dx\le  0 \qquad \text{for every $c>0$,}
$$
and by (ii) we conclude that $w = u+ \gamma(s_0) \in \cM$.\\
(iv) Let $u \in W^{1,2}(\Omega)$. The uniqueness of $c=c(u) \in \R$ with $u+c
\in \cM$ is an immediate consequence of (i). To see the existence, we
note that $\sgn(u\pm c) \to \pm 1$ as $c \to +\infty$ a.e. in
$\Omega$. Hence, by Lebesgue's theorem, there exists $c_0>0$ with 
$$
\pm \int_{\Omega} \sgn(u\pm c_0)\,dx > 0. 
$$
Applying (iii) to the path $s \mapsto (1-2s)c_0$ now yields the
existence of $c \in [-c_0,c_0]$ such that $u + c \in \cM$. The
continuity follows similarly as (but more easily than) in
Lemma~\ref{sec:neum-bound-cond-3}.\\
(v) Suppose by contradiction that $c(u)<c(v)=:c$. We then have 
$u \le u+c \le v+c$ in $\Omega$ and therefore
$$
\int_{\Omega}\sgn_-(u+c)\,dx \le \int_{\Omega}\sgn_-(v+c)\,dx \le 0
$$
and 
$$
\int_{\Omega}\sgn_+(u+c)\,dx \ge \int_{\Omega}\sgn_+(u)\,dx \ge 0.
$$
By (\ref{eq:29}), we then have $u+c \in \cM$, with contradicts the
uniqueness statement in (iv).
\end{proof}

We now consider the variational problem related to the minimax value 
\[ 
m:= \inf_{\cM} \varphi = \inf_{u \in W^{1,2}(\Omega)} \sup_{c \in \real} \varphi (u+c) 
\]
Note that the second equality is an immediate consequence of
Lemma~\ref{sec:case-q=1}(i), (iv). Note also that $m<0$, since for
every $u \in \cM \setminus \{0\}$ we have 
\begin{equation}
  \label{eq:26}
\phi(t^*(u)u) = -\frac{1}{2} \frac{\|u\|_{L^1(\Omega)}^2}{\int_{\Omega} |\nabla
u|^2\,dx} <0 \qquad \text{with $t^*(u)=
\frac{\|u\|_{L^1(\Omega)}^2}{\int_{\Omega}|\nabla u|^2\,dx}$,}
\end{equation}
whereas $t^*(u)u \in \cM$.

The main result of this section is the following. 

\begin{thm}
\label{sec:case-q=1-1}
The value $m<0$ is attained by $\phi$ on $\cM$. Moreover, every $u \in
\cM$ with $\phi(u)=m$ is a nontrivial solution of \eqref{problem_q=1}
such that its zero set $\{u=0\}  \subset \Omega$ has vanishing
Lebesgue measure. 
\end{thm}

The proof of this Theorem is split in two steps. We first show the following.

\begin{lem}
\label{sec:neum-bound-cond-1-1-1}
The functional $\varphi$ attains the value $m<0$ on $\cM$. Moreover,
every minimizer $u$ of $\varphi$ on $\cM$ is a sign changing solution
of \eqref{problem_q=1}.
\end{lem}

\begin{proof}
We first note that for all $u \in \cM$ we have, by Lemma~\ref{sec:case-q=1}(i)
$$
\|u\|_{L^1(\Omega)}^2 = \min_{c \in \R} \|u+c\|_{L^1(\Omega)}^2 \le 
|\Omega| \min_{c \in \R} \|u+c\|_{L^2(\Omega)}^2 \le  |\Omega| \mu^{-1}_2 \int_{\Omega}|\nabla u|^2\,dx,
$$
where $\mu_2>0$ is the first nontrivial Neumann eigenvalue of
$-\Delta$ on $\Omega$. As a consequence, the functional $\varphi$ is coercive on $\cM$. Let
$(u_n)_n  \subset \cM$ be a minimizing sequence for $\varphi$. Then
$(u_n)$ is bounded, and we may pass to a subsequence such that $u_n
\weak \tilde u \in  W^{1,2}(\Omega)$. Then $u_n \to \tilde u$ in $L^1(\Omega)$,
$$
\int_{\Omega}|\nabla \tilde u|^2\,dx \,{\le}\, \liminf_{n \to \infty} \int_{\Omega}|\nabla u_n|^2\,dx 
$$
and 
$$
\int_{\Omega}|\tilde u + c| \,dx = \lim_{n \to \infty} \int_{\Omega}|u_n +
c|\,dx \le \lim_{n \to \infty} \int_{\Omega}|u_n|\,dx = \int_{\Omega}
|\tilde u|\,dx \qquad \text{for all $c \in \R$.}
$$
Consequently, $\tilde u \in \cM$ by Lemma~\ref{sec:case-q=1}(i) and
(iv), and 
$$
\phi(\tilde u) \le \liminf_{n \to \infty} \phi(u_n) = m.
$$
By definition of $m$, equality holds and thus $\phi$ attains its minimum on $\cM$.\\
Next, we let $u \in \cM$ be an arbitrary minimizer for $\varphi$ on
$\cM$, and we show that $u$ is a solution of \eqref{problem_q=1}.
We first show that 
\begin{equation}
  \label{eq:22}
\int_{\Omega}\sgn_-(u) v\,dx \le  \int_{\Omega}\nabla u \nabla v\,dx \le \int_{\Omega}\sgn_+(u) v\,dx
\qquad \text{for every $v \in W^{1,2}(\Omega)$, $v \ge 0$}.  
\end{equation}
Arguing by contradiction, we assume that there exists $v \in
W^{1,2}(\Omega)$, $v \ge 0$ such that 
\begin{equation}
  \label{eq:28}
\int_{\Omega} \sgn_-(u) v\,dx >  \int_{\Omega} \nabla u \nabla v\,dx.
\end{equation}
Note that for $a,c\in \R$ we have $|a+c| \ge |a| + \sgn_-(a)c$. Hence for every $c \le 0$, $t \ge 0$ we have 
$$
\|u+c+tv\|_{L^1(\Omega)} \ge \|u\|_{L^1(\Omega)} + \int_{\Omega}\sgn_-
(u)(c+tv)\,dx \ge \|u\|_{L^1(\Omega)} + t \int_{\Omega}\sgn_-
(u)v\,dx, 
$$
where the last inequality follows from the fact that $\int_{\Omega}
\sgn_-(u)\,dx \le 0$ since $u \in \cM$. Since $v \ge 0$ implies that $c(u+tv) \le 0$ for $t >
0$ by Lemma~\ref{sec:case-q=1}(v), we find that 
\begin{align*}
\varphi(u+tv+c(u+tv)) &=  \frac{1}{2}\int_{\Omega}|\nabla (u+tv)|^2\,dx - 
\|u+tv +c(u+tv)\|_{L^1(\Omega)}\\ 
&\le \frac{1}{2}\int_{\Omega}|\nabla u|^2\,dx - \|u\|_{L^1(\Omega)} +t \Bigl(\int_{\Omega} \nabla u \nabla v\,dx - \int_{\Omega}\sgn_- (u)
v\,dx\Bigr)+ \frac{t^2}{2}\int_{\Omega}|\nabla v|^2\,dx\\ 
&=m +t \Bigl(\int_{\Omega} \nabla u \nabla v\,dx - \int_{\Omega}\sgn_- (u)
v\,dx\Bigr)+ \frac{t^2}{2} \int_{\Omega}|\nabla v|^2\,dx\\
&<m \qquad \text{for $t>0$ sufficiently small}
\end{align*}
by (\ref{eq:28}). This contradicts the definition of $m$. Hence the first inequality in
(\ref{eq:22}) holds, and the second inequality is proved by a similar
argument. As a consequence, we have
$$
\int_{\Omega}\nabla u \nabla v \,dx \le \int_{\Omega} \bigl[\sgn_+(u)v^+ -
\sgn_-(u)v^-\bigr]\,dx \le \|v\|_{L^1(\Omega)} \qquad \text{for every $v \in W^{1,2}(\Omega)$.} 
$$
Consequently, the distributional Laplacian $\Delta u:
C_0^\infty(\Omega) \to \R$ is continuous with respect to the
$L^1(\Omega)$-norm and is therefore represented by
a function $- w \in L^\infty(\Omega)$ satisfying 
\begin{equation}
  \label{eq:23}
\sgn_-(u) \le w \le \sgn_+(u).  
\end{equation}
Then, by elliptic regularity theory, it follows that $u \in W^{2,p}_{loc}(\Omega)$ for
all $p \in (1,\infty)$ with $-\Delta u =w$. Moreover, we have $\nabla u \equiv 0$ and
$w=-\Delta u \equiv 0$ a.e. on the set $\{u=0\}$ (see
e.g.  \cite[Lemma 7.7]{gt}). Hence, by
(\ref{eq:23}) we may assume that $w=\sgn(u)$, and thus $u$ is a
solution of \eqref{problem_q=1}. Finally, to show that $u$ is sign
changing, we first note that $u
\not = 0$ since $\phi(u)=m<0$. Suppose by contradiction that $u \ge
0$, then $u$ is also superharmonic by \eqref{problem_q=1}, and hence
$u>0$ in $\Omega$ by the strong maximum principle, which contradicts
the fact that $u \in \cM$. Similarly, we get a contradiction assuming
that $u \le 0$. Hence $u$ changes sign in $\Omega$.
\end{proof}

In order to complete the proof of Theorem~\ref{sec:case-q=1-1}, we
need to show that every minimizer $u \in \cM$ of $\phi$ has the
unique continuation property, that is, $\{u=0\}$ has measure zero. The
argument is similar as in the case $q>1$, but changes are
required at some points. We start with the following. 

\begin{prop}
\label{sec:some-estim-relat-2-q=1}
Let $u$ be a solution of \eqref{problem_q=1}, and suppose that $x_0 \in \Omega$ is a point of density one for the set
$u^{-1}(0) \subset \Omega$. Then $u(x) = o(|x|^{2})$
as $x \to x_0$.   
\end{prop}

\begin{proof}
The argument is similar as the proof of
Proposition~\ref{sec:some-estim-relat-2}. Without loss, we assume
that $x_0=0$, and we extend $u$ to all of $\R^N$ by setting $u \equiv
0$ on $\R^N \setminus \Omega$. Applying Lemma~\ref{sec:some-estim-relat-1} as in the proof of
Proposition~\ref{sec:some-estim-relat-2} with $\alpha=1$, $\beta=2$,
we see that the function 
$$
g:(0,\infty) \to \R,\qquad g(r):= r^{-2}\sup \limits_{|x|=r}|u(x)|
$$
is bounded. To show that $g(r) \to 0$ as $r \to
0$, we argue by contradiction, assuming that there exists $\eps>0$ and a sequence of radii $r_n
>0$, $n \in \N$ such that $r_n \to 0$ as $n \to \infty$ and $g(r_n)
\ge \eps$ for all $n \in \N$. We may assume that $B_{2r_n}(0) \subset
\Omega$ for all $n \in \N$. We then consider $\Omega_0:= B_2(0)
\setminus \overline{B_{\frac{1}{2}}(0)}$ and the functions $w_n: \Omega_0 \to \R$,
$w_n(x)= r_n^{-2}
u(r_n x)$ for $n \in \N$ which are uniformly bounded in $\Omega_0$.
For $n \in \N$, $w_n$ solves 
$$
-\Delta w_n = \sgn(w_n) \qquad \text{in $B_2(0) \setminus
  B_{\frac{1}{2}}(0)$},
$$
and there exists a sequence of points $x_n \in S^1$, $n
\in \N$ such that $w(x_n)=g(r_n) \ge \eps$ for all $n \in \N$. Using elliptic regularity theory again, we may pass to a subsequence such that 
\begin{align*}
&x_n \to \bar x \in S^1;\\
&\text{$w_n \to w$ in $C^1_{loc}(\Omega_0)$,}\\
&\sgn(w_n) \weak^* f \qquad \text{in $L^\infty(\Omega_0)=(L^1(\Omega_0))^*$,}
\end{align*}
where $w \in C^{1,\alpha}_{loc}(\Omega_0)$
is a weak solution of $-\Delta w = f$ in $\Omega_0$ such that
$w(\bar x) \ge \eps$. In particular, $w \not \equiv 0$. 
However, by the same argument as in the proof of Proposition~\ref{sec:some-estim-relat-2}, the sets $A_n:= \{x \in B_2(0) \setminus B_{\frac{1}{2}}(0)\::\:
w_n(x) \not= 0\}$ satisfy $\lim \limits_{n \to \infty}|A_n|=0$ and
therefore $w \equiv 0$ a.e. in $B_2(0) \setminus
B_{\frac{1}{2}}(0)$. This yields a contradiction, and thus we conclude
that $g(r) \to 0$ as $r \to 0$, as required.
\end{proof}

Next, for $r>0$ and $q=1$, we consider the functional
$\phi_r: W^{1,2}(B_r(0)) \to \R$ and the rescaling map $T_r:
W^{1,2}(B_1(0)) \to W^{1,2}(B_r(0))$ defined in (\ref{eq:24}),
(\ref{eq:25}), respectively.
We then have the following. 

\begin{prop}
\label{sec:unique-continuation-1-q=1}
Let $u \in W^{1,2}(\Omega)$ be a solution of \eqref{problem_q=1} in $\Omega$,
and let $x_0 \in \Omega$ be a point of density one for $u^{-1}(0) \subset
\Omega$. Moreover, let $K \subset W^{2,2}_0(B_1(0))$ be a compact set
such that 
$$
c_K:= \sup_{v \in K}\phi_1(v)<0,
$$
Then there exists $r_0>0$ with the following property:\\
$B_{r_0}(x_0)$ is contained in $\Omega$, and for every $r \in (0,r_0)$ and every $v \in K$ we have $\phi(u+v_r)<\phi(u)$,
where $v_r \in W^{2,2}_0(B_r(x_0)) \subset W^{2,2}_0(\Omega)$ is defined by 
$$
v_r(x)= \left\{
  \begin{aligned}
   &[T_r v](x-x_0)&&\qquad x \in B_{r}(x_0);\\
   &0 &&\qquad x \in \Omega \setminus B_r(x_0).     
  \end{aligned}
\right.
$$
\end{prop}

The proof is somewhat different than the proof of
Proposition~\ref{sec:unique-continuation-1}. Note that we need the
stronger assumption $K \subset W^{2,2}_0(B_1(0))$ here. This
assumption is not optimal but suffices for our purposes. 

\begin{proof}
 Without loss, we may assume that $x_0=0 \in \Omega$. Let $v \in K$. We
then have 
\begin{align}
\nonumber
\phi(u+v_r)&=\phi(u)+\phi_r(v_r)+ \int_{B_r(0)}\nabla u \nabla v_r\,dx
- \int_{B_r(0)}\Bigl(|u+v_r|-|u|-|v_r|\Bigr))\,dx\\
&\le \phi(u)+\phi_r(v_r)+ \int_{B_r(0)}\nabla u \nabla v_r\,dx +
 2 \int_{B_r(0)}|u|\,dx .  \label{eq:8-1}
\end{align}
Note that, since $v_r \in W^{2,2}_0(B_r(0))$, we have 
\begin{align*}
\Bigl|\int_{B_r(0)}\nabla u \nabla v_r\,dx\Bigr| &= \Bigl|
\int_{B_r(0)}u \Delta v_r\,dx\Bigr| \le r^{N}
\Bigl| \int_{B_1(0)}u(rx) [\Delta v](x) \,dx\Bigr| \\
&= o(r^{N+2})\|\Delta v\|_{L^1(B_1(0))}= o(r^{N+2})\qquad \text{as $r
  \to 0$.}   
\end{align*}
Moreover,
$$
\int_{B_r(0)}|u|\,dx = o(r^{2})|B_r(0)|= o(r^{N+2})\qquad \text{as $r
  \to 0$.}
$$
Finally, since $\nabla v_r(x)= r  \nabla v(\frac{x}{r})$ for $x
\in B_r(0)$, we find
that   
$$
\phi_r(v_r)= r^{2} \Bigl(\frac{1}{2} \int_{B_r(0)}\Bigl|\nabla
v\left(\frac{x}{r}\right)\Bigr|^2\,dx
- \int_{B_r(0)}\Bigl|v\left(\frac{x}{r}\right)\Bigr|\,dx \Bigr)=r^{N+2}\phi_1(v).
$$
Inserting these estimates in (\ref{eq:8-1}), we obtain 
$$
\phi(u+v_r)\le \phi(u)+ r^{N+ 2}\Bigl(\phi_1(v)+o(1)\Bigr) \le 
\phi(u)+ r^{N+ 2}\Bigl(c_K+o(1)\Bigr) .
$$
Since $K$ is compact, it is easy to see that these estimates are
uniform in $v \in K$. Since moreover $c_K<0$, the claim follows.
\end{proof}

\begin{thm}
\label{sec:neum-bound-cond-4-q=1}
Let $u$ be a minimizer of $\phi$ on $\cM$. Then $u^{-1}(0)$ has vanishing Lebesgue
measure.   
\end{thm}

\begin{proof}
Suppose by contradiction that $|u^{-1}(0)|>0$. Then there exists a
point $x_0$ of density one for the set $u^{-1}(0)$. Without loss of generality, we can suppose $x_0=0$. We fix two arbitrary nonnegative nontrivial functions $v_1,v_2 \in
C_c^2(B_1(0))$ with disjoint support, and consider the path 
$$
\gamma:[0,1] \to C_c^2(B_1(0)),\qquad \gamma(s)= t^{*}((1-s) v_1 -
s v_2)\cdot ((1-s) v_1 -s v_2),
$$
where the map $t^*$ is defined in (\ref{eq:26}). It is clear
that $\varphi(\gamma(s))<0$ for $s \in [0,1]$. We also define
$$
\gamma_r := T_r \circ \gamma :[0,1] \to C_c^2(B_r(0)).
$$
Applying Proposition~\ref{sec:unique-continuation-1-q=1} to the compact set $K:=
\gamma([0,1]) \subset W^{2,2}_0(B_1(0))$, we may fix $r>0$
sufficiently small such that $\phi(u+ \gamma_r(s))<\phi(u)$ for every
$s \in [0,1]$. Moreover, by making $r$ smaller if necessary and using again the fact
that $x_0=0$ is a point of density one for the set $u^{-1}(0)=0$, we
may assume that 
$$
|\{u=0\} \cap \{\gamma_r(0) > 0\}|>0   \qquad \text{and}\qquad 
|\{u=0\} \cap \{\gamma_r(1) < 0\}|>0.
$$
As a consequence of these
inequalities and the fact that $u$ is a
solution of \eqref{problem_q=1}, we find that 
$$
\int_{\Omega}\text{sgn}(u+\gamma_r(0))\,dx \:>\: 
\int_{\Omega}\text{sgn}(u)\,dx=0 \:>\: \int_{\Omega}\text{sgn}(u+\gamma_r(1))\,dx <0.
$$
By Lemma~\ref{sec:case-q=1}(iii), there exists $s_0 \in [0,1]$ such that
$u+\gamma_r(s_0) \in \cM$. 
This however contradicts the assumption that $u$ is a minimizer of $\phi$ in $\cM$.
\end{proof}

In the following, we restrict our attention to the case where $\Omega$
is a radial bounded domain in $\R^N$. In this case we also consider
$$
\cM_r:= \{u \in \cM\::\: \text{$u$ radial}\}\qquad \text{and}\qquad m_r:=
\inf_{\cM_r} \phi,
$$
so that $\cM_r \subset \cM$ and $m_r \ge m$.

\begin{thm}
\label{sec:case-q=1-3}
Let $\Omega \subset \R^N$ be a radial bounded domain. Then we have:
\begin{itemize}
\item[(i)] If $u \in \cM$ satisfies $\phi(u)=m$, then $u$ is foliated
  Schwarz symmetric. 
\item[(ii)] If $u \in \cM_r$ satisfies $\phi(u_r)=m_r$, then $u$ is strictly
monotone in the radial variable.
\item[(iii)] If $\Omega =B$ is the unit ball in $\R^N$, we have 
\begin{equation}
  \label{eq:30bis}
m \le  -\frac{\pi}{18} <m_r = -\pi \left( -\frac{1}{16} + \frac{1}{8}\ln{2} \right),
\end{equation}
if $N=2$ and, for $N \geq 3$,
\begin{equation}
  \label{eq:30}
m \le  -\omega_N \frac{N-2}{2N^2(N-1)} <m_r = -\frac{\omega_N}{2} \frac{(2^{-\frac{2}{N}}-1)N+2^{1-\frac{2}{N}}}{(N-2)(N+2)},
\end{equation}
where, as usual, $\omega_N$ denotes the measure of $|B|$. Hence every minimizer $u \in \cM$ of $\phi|_{\cM}$ is a nonradial function.
\end{itemize}
\end{thm}

\begin{proof}
(i) The proof is very similar to the proof of Theorem~\ref{sec:neum-bound-cond-5}, and we use the notation introduced
there. Pick $x_0 \in
\Omega \setminus \{0\}$ with $u(x_0)=\max \{u(x)\::\: |x|=|x_0|\}$.
We put $p:=\frac{x_0}{|x_0|}$, and we let $H \in \cH_p$. As in the proof of
Theorem~\ref{sec:neum-bound-cond-5}, it suffices to show that $u
\equiv u_H$ on $H \cap \Omega$. Since
$$
\int_{\Omega} |\nabla u_H|^2\:dx = \int_{\Omega} |\nabla
u|^2\:dx,\quad
\int_{\Omega} |u_H|\:dx = \int_{\Omega} |u|\:dx \quad
\text{and}\quad \int_{\Omega} \sgn_{\pm}(u_h)\:dx = \int_{\Omega}\sgn_{\pm}(u) \:dx,
$$
we find that $u_H \in \cM$ and $\phi(u_H)=\phi(u)$, so that $u_H$ is
also a minimizer of $\phi$ on $\cM$. Hence, by Theorem~\ref{sec:neum-bound-cond-1}, both $u$ and $u_H$ are
solutions of \eqref{problem_q=1}. Therefore $w:=u_H-u$ is a
nonnegative function in $\Omega \cap H$ satisfying 
$$
-\Delta w= \sgn(u_H) -\sgn(u) \ge 0\qquad \text{in
    $H \cap \Omega$.}
$$
The strong maximum principle then implies that either $w \equiv 0$ or $w>0$ in $H \cap \Omega$. The
latter case is ruled out since $x_0 \in H \cap \Omega$ and
$w(x_0)=u_H(x_0)-u(x_0)=0$ by the choice of $x_0$. We therefore obtain
$w \equiv 0$ and hence $u \equiv u_H$ on $H \cap \Omega$, as required.\\
(ii) We only consider the case where $\Omega=B$ is the unit ball in $\R^N$;
the proof in the case of an annulus is similar. Let $u \in \cM_r$
satisfy $\phi(u_r)=m_r$. Then, by similar arguments as in the proof of
Theorem~\ref{sec:neum-bound-cond-1-1}, $u$ is a radial 
solution of \eqref{problem_q=1}. Suppose by contradiction that there
exists $r_0 \in (0,1)$ such that $u'(r_0)=0$. We claim that we can
choose $r_0$ minimally, i.e., such that  
\begin{equation}
  \label{eq:39}
u'(r) \not=0 \quad \text{for $r \in (0,r_0)$.}  
\end{equation}
Indeed, suppose by contradiction that there exists a sequence of $r_n
\in (0,1)$, $n \in \N$ such that $u'(r_n)=0$ for all $n$ and $r_n \to
0$ as $n \to \infty$. Without loss, we may assume that $r_n >r_{n+1}$
for every $n$. Since the function $r \mapsto r^{N-1}u'(r)$ is
strictly monotone on every interval on which $u$ has no zero, we
conclude that there exists $s_n \in (r_{n+1},r_n)$ such that
$u(s_n)=0$. Since $u$ is of class $\cC^1$, we therefore conclude that
$u(0)=u'(0)=0$. This however implies that $u \equiv 0$, since the absolutely continuous function $r \mapsto h(r)=
\frac{u'(r)^2}{2} +|u(r)|$ is decreasing on $[0,1]$, which
follows from the fact that 
\begin{equation}
  \label{eq:35}
h'(r)= u'(r)[u''(r)+\sgn(u(r))] =-\frac{(N-1)}{r}u'(r)^2 \le 0 \qquad
\text{for a.e. $r \in (0,1)$.}
\end{equation}
We thus conclude that we can choose $r_0 \in (0,1)$ such that
(\ref{eq:39}) holds. It is then easy to see (e.g. by using the Hopf
boundary lemma) that $u(r_0)\not=0$. Let $\Omega_1:=
B_{r_0}(0)$ and $\Omega_2:= B \setminus B_{r_0}(0)$. Then $u$ solves
\eqref{problem_q=1} both on $\Omega_1$ and $\Omega_2$, so that 
$$
\int_{\Omega_1} \sgn(u)\,dx = 0 = \int_{\Omega_2}\sgn(u)\,dx.
$$
We now define 
$$
v \in W^{1,2}(B),\qquad v(x) = \left \{
  \begin{aligned}
  &2u(r_0)-u(x),&&\qquad x \in \Omega_1;\\
  &u(x),&&\qquad x \in \Omega_2.   
  \end{aligned}
\right.
$$
Moreover, we let $c=c(v) \in \R$ be given by
Lemma~\ref{sec:case-q=1}(iv), so that $w:=v+c \in \cM_r$. 
We then have 
\begin{equation}
  \label{eq:43}
\int_{B}|\nabla w|^2\,dx = \int_{B}|\nabla v|^2\,dx = \int_{B}|\nabla u|^2\,dx.
\end{equation}
Moreover,
\begin{equation}
  \label{eq:40}
\int_{\Omega_1}|w|\,dx = \int_{\Omega_1} |u-2u(r_0)-c|\,dx \ge 
\int_{\Omega_1} |u|\,dx
\end{equation}
by Lemma~\ref{sec:case-q=1}(i) applied to the domain
$\Omega_1$. Similarly, 
\begin{equation}
  \label{eq:41}
\int_{\Omega_2}|w|\,dx = \int_{\Omega_2} |u+c|\,dx \ge
\int_{\Omega_2} |u|\,dx
\end{equation}
by Lemma~\ref{sec:case-q=1}(i) applied to the domain
$\Omega_2$. Moreover, if equality holds in both (\ref{eq:40}) and
(\ref{eq:41}), then Lemma~\ref{sec:case-q=1}(i) implies that 
$2u(r_0)+c=0$ and $c=0$, hence $u(r_0)=0$ contrary to what we have
seen earlier. Hence at least one of the inequalities (\ref{eq:40}) and
(\ref{eq:41}) is strict, so that 
\begin{equation}
  \label{eq:42}
\int_{\Omega}|w|\,dx> \int_{\Omega}|u|\,dx  
\end{equation}
and therefore $\phi(w)<\phi(u)=m_r$ as a consequence of (\ref{eq:43}) and
(\ref{eq:42}). This contradicts the definition of $m_r$, and thus the
proof is finished.\\
(iii) Similar arguments as in the proof of Lemma~\ref{sec:neum-bound-cond-1-1-1} show that $m_r$ is attained by a
radial solution $u \in \cM_r$ of \eqref{problem_q=1}. By (ii), we know
that $u$ is strictly monotone in the radial variable, and therefore it
is up to sign uniquely given by
\[ u(r) = \left\{ \begin{array}{c l} \displaystyle -\frac{1}{4}r^2 + \frac{1}{8}  & \text{for } 0 \leq r \leq a := \frac{1}{\sqrt{2}} \\ \displaystyle
                   -\frac{1}{2}\ln{r} + \frac{1}{4}r^2 - \frac{1}{8} - \frac{1}{4}\ln{2}  & \text{for } a \leq r \leq 1 
                  \end{array} \right.
\]
for $N=2$, and by
$$ 
u(r) = \left\{ \begin{array}{c l} \displaystyle \frac{1}{2N}(a^2 - r^2) & \text{for } 0 \leq r \leq a := 2^{-\frac{1}{N}} \\ \displaystyle
                   \frac{1}{2N} \left( \frac{2}{N-2}r^{2-N} + r^2 - \frac{2^{-\frac{2}{N}}(N+2)}{N-2} \right) & \text{for } a \leq r \leq 1 
                  \end{array} \right.
$$
for $N \geq 3$. Note here that the value $a$ is determined by the
condition that $|\{u>0\}|=|\{u<0\}|$. For $N=2$, we have
$$
\int_B |\nabla u|^2 =  \int_B |u| =  2\pi \left( -\frac{1}{16} + \frac{1}{8}\ln{2} \right)
$$
and therefore
$$
m_r=\phi(u)= -\pi \left( -\frac{1}{16} + \frac{1}{8}\ln{2} \right),
$$
thus showing the right inequality in \eqref{eq:30bis}. For $N\geq 3$,
$$
\int_B |\nabla u|^2 =  \int_B |u| =  \frac{\omega_N}{2} \left(
  \frac{1}{N-2}2^{\frac{N-2}{N}} - \frac{1}{N+2} -
  \frac{1}{N-2}\right) = \omega_N \frac{(2^{-\frac{2}{N}}-1)N+2^{1-\frac{2}{N}}}{(N-2)(N+2)}
$$
and therefore 
$$
m_r=\phi(u)= -\frac{\omega_N}{2} \frac{(2^{-\frac{2}{N}}-1)N+2^{1-\frac{2}{N}}}{(N-2)(N+2)},
$$
which is the right equality in (\ref{eq:30}). To see the left
inequalities in \eqref{eq:30bis} and (\ref{eq:30}), we consider the functions $x \mapsto
u_s(x)= x_1|x|^s$ for $s >-\frac{N}{2}$. We then have 
$$
\int_{B} |u_s|\,dx \ge \int_{B}x_1^2 |x|^{s-1}\,dx = \frac{1}{N}
\int_{B}|x|^{s+1}\,dx = \frac{\omega_N}{N+s+1}
$$
and 
$$
\int_{B} |\nabla u_s|^2 \,dx=
\int_{B} \Bigl(|x|^{2s}+(s^2+2s)x_1^2|x|^{2(s-1)} 
\Bigr)\,dx=\omega_N\frac{N+s^2+2s}{N+2s}.
$$
We now have $v_s \in \cM$ for the function $v_s= c_s u_s$ with 
$c_s:= \frac{\int_{B}|u_s|\,dx}{\int_{B}|\nabla u_s|^2\,dx}$ and
therefore 
$$
m \le \phi(v_s)= -\frac{1}{2} \frac{\Bigl(\int_{B}|u_s|\,dx\Bigr)^2}{\int_{B} |\nabla u_s|^2 \,dx} \leq -
\omega_N \frac{N+2s}{2(N+s^2+2s)(N+s+1)^2}.
$$
Thus the left inequalities in \eqref{eq:30bis} and (\ref{eq:30}) follow by choosing $s=0$ if $N=2$, and $s=-1$ if $N \geq 3$.\\ 
Finally, for the middle inequality in (\ref{eq:30}) we need to prove: 
\begin{equation} \label{equationwiththeN}
\frac{(2^{-\frac{2}{N}}-1)N+2^{1-\frac{2}{N}}}{(N-2)(N+2)}<\frac{N-2}{N^2(N-1)}
\end{equation}
for $N \geq 3$, which is equivalent to
\begin{align*} & (2^{-\frac{2}{N}}-1)N^4 + 2^{-\frac{2}{N}}N^3 + 2(1-2^{-\frac{2}{N}})N^2 + 4N - 8 < 0 \\ & \Leftrightarrow N^3 \left[(2^{-\frac{2}{N}}-1)N + 2^{-\frac{2}{N}} + \frac{2}{N}(1-2^{-\frac{2}{N}})\right] + 4N -8 < 0. \end{align*}
To this aim, we will prove that the function $h:[3, +\infty) \to \real$ defined by
\[ h(t) := h_1(t) + h_2(t) + h_3(t) = (2^{-\frac{2}{t}}-1)t + 2^{-\frac{2}{t}} + \frac{2}{t}(1-2^{-\frac{2}{t}}) \]
is maximized for $t=3$, and that $h(3)<-1$. Inequality \eqref{equationwiththeN} will be implied by $-N^3 + 4N - 8 < 0$, which is true for $N \geq 3$.
The function $h_1$ is such that
\[ h_1'(t) = \frac{4^{-\frac{1}{t}}\left(-4^{\frac{1}{t}}t + t + \ln{4} \right)}{t}, \qquad h_1''(t) = \frac{4^{\frac{t-1}{t}}(\ln{2})^2}{t^3}. \]
Since $h_1'(t)\to 0$ as $t \to +\infty$ and $h_1''(t) > 0$ for $t>0$, $h_1$ is monotone decreasing. $h_2$ is clearly strictly decreasing, while $h_3$ satisfies
\[ h_3'(t) = -\frac{2^{\frac{t-2}{t}}\left((4^{\frac{1}{t}}-1)t + \ln{4} \right)}{t^3} < 0 \]
for $t>0$. It is then easily verified that $h(3) = \frac{1}{3}(5 \sqrt[3]{2} - 7) < 1$.
\end{proof}

\bibliographystyle{amsplain}

\end{document}